\documentclass[onefignum,onetabnum]{siamart251216}

\usepackage{lipsum}
\usepackage{amsfonts}
\usepackage{graphicx}
\usepackage{epstopdf}
\usepackage{algorithmic}
\ifpdf
  \DeclareGraphicsExtensions{.eps,.pdf,.png,.jpg}
\else
  \DeclareGraphicsExtensions{.eps}
\fi

\newsiamremark{remark}{Remark}
\newsiamremark{hypothesis}{Hypothesis}
\crefname{hypothesis}{Hypothesis}{Hypotheses}
\newsiamthm{claim}{Claim}

\headers{A General Framework for Preconditioning \NoCaseChange{$f(A)b$}}{G. Ramirez-Hidalgo}

\title{A General Framework for Preconditioning \NoCaseChange{$f(A)b$}\thanks{Submitted to the editors DATE.
\funding{This work was funded by the Fog Research Institute under contract no.~FRI-454.}}}

\author{Dianne Doe\thanks{Imagination Corp., Chicago, IL 
  (\email{ddoe@imag.com}, \url{http://www.imag.com/\string~ddoe/}).}
\and Paul T. Frank\thanks{Department of Applied Mathematics, Fictional University, Boise, ID 
  (\email{ptfrank@fictional.edu}, \email{jesmith@fictional.edu}).}
\and Jane E. Smith\footnotemark[3]}

\usepackage{amsopn}

\makeatletter
\newcommand*{\addFileDependency}[1]{
  \typeout{(#1)}
  \@addtofilelist{#1}
  \IfFileExists{#1}{}{\typeout{No file #1.}}
}
\makeatother

\usepackage{amsmath}
\usepackage{amssymb}

\usepackage{xcolor}
\usepackage{graphicx}
\usepackage{subcaption} 

\newtheorem{assumption}[theorem]{Assumption}

\begin{document}


\title{An Integral-Based Framework for Preconditioning \lowercase{$f$}$(A)$\lowercase{$b$}}

\author{Gustavo Ramirez-Hidalgo\thanks{Jülich Supercomputing Centre, Forschungszentrum Jülich GmbH, Wilhelm-Johnen-Straße 52428, Jülich, Germany (\email{g.ramirez.hidalgo@fz-juelich.de}).}}

\headers{An Integral-Based Framework for Preconditioning \lowercase{$f$}$(A)$\lowercase{$b$}}{G. Ramirez-Hidalgo}

\maketitle

\begin{abstract}
    The computation of the action of a matrix function on a vector, $f(A)b$, is a major computational bottleneck for large, sparse matrices, particularly when unfavorable spectral distributions cause standard Krylov subspace methods to stagnate. In this work, we propose a unified framework for preconditioning $f(A)b$ based on the Cauchy integral representation of the matrix function. By exploiting shift-invariance properties, we decouple the preconditioner evaluation from the Krylov subspace generation. We develop this framework in two distinct directions. First, for rational shift-and-invert preconditioning, we resolve a fundamental trade-off between optimal spectral compression and finite-precision instability. We achieve this by formulating a closed-form extraction stabilized via Double Modified Gram-Schmidt reorthogonalization, which eliminates the formation of spurious phantom poles. Second, we present a matrix-free polynomial approach. To ensure numerical stability, we isolate the continuous numerical quadrature step using a Schur decomposition of the projected Hessenberg matrix. To further stabilize the integration near contour singularities and accelerate overall convergence, we incorporate an exact LR-deflation scheme targeting the critical low modes of the preconditioned operator. We analyze the asymptotic stability and proximity to singularity of these methods, and present numerical experiments demonstrating their efficiency on the 2D Laplacian with $f=\textrm{exp}$, and a highly ill-conditioned Wilson-Dirac operator from lattice quantum chromodynamics with $f=\textrm{sign}$, although the framework can be in principle used with any $f$ and it is particularly beneficial when applying $f(A)b_{i}$ with many different vectors $b_{i}$.
\end{abstract}



\section{Introduction}

The computation of the action of a matrix function on a vector, denoted as $f(A)b$, is a ubiquitous task in scientific computing. It appears centrally in applications ranging from the solution of time-dependent partial differential equations (PDEs) via exponential integrators to lattice quantum chromodynamics (QCD) and complex network analysis. When the matrix $A \in \mathbb{C}^{N \times N}$ is large and sparse, direct computation of the full dense matrix $f(A)$ is computationally intractable. Consequently, standard iterative approaches, such as the Arnoldi and Lanczos methods, construct a Krylov subspace approximation of the solution. However, for functions possessing singularities near the spectrum of $A$, or for matrices exhibiting highly unfavorable spectral distributions---such as ill-conditioned systems arising from fine PDE discretizations---the convergence of standard Krylov methods can be unacceptably slow.

In the realm of standard linear systems ($Ax = b$), preconditioning is a mature and indispensable tool for accelerating convergence. Conversely, preconditioning for general matrix functions remains a challenging and relatively sparse frontier. The fundamental difficulty lies in the fact that, unlike a linear system, $f(A)b$ is not defined by a simple algebraic residual equation that can be easily transformed by an arbitrary operator.

The primary objective of this work is to start closing this gap by establishing a rigorous, highly scalable framework for preconditioning $f(A)b$. We accomplish this by transforming the spectral properties of the underlying operator using shifted preconditioners tailored directly to the Cauchy integral representation of the matrix function. 

A key innovation of our framework is the reliance on shift-invariance properties across the contour integration, which decouples the preconditioner evaluation from the Krylov subspace generation. We explore two distinct branches of this framework:
\begin{enumerate}
    \item \textit{Shift-and-Invert Preconditioning:} We prove that utilizing a rational preconditioner yields a mapped integral representation. While standard algebraic extractions suffer from finite-precision noise amplification due to spectral clustering, we formulate a closed-form approach stabilized by Double Modified Gram-Schmidt (DMGS) reorthogonalization that shields the approximation from finite-precision phantom poles, and with a cap on the shift of the shift-and-invert preconditioner to avoid singularities in some small-subspace dense inversions involved in the method.
    \item \textit{Polynomial Preconditioning:} To circumvent the possibility of the shift-and-invert preconditioner being overkill, we formulate a matrix-free polynomial approach. The contour-based framework presented here allows us to bypass the possibly problematic algebraic sign ambiguities that can arise in other polynomial preconditioning approaches when computing functions like the inverse square root. We introduce an algorithmic design utilizing Schur decompositions and dual-regime scalings to insulate the numerical quadrature from floating-point instabilities. Furthermore, to stabilize polynomial integration near singular contour regions, we integrate an exact LR-deflation scheme for the critical low modes of the preconditioned operator, which at the same time accelerates the shifted linear systems induced by the contour integration, leading to an overall reduction in computational cost.
\end{enumerate}

We provide, in Sections \ref{sect:Arnoldi_method_for_fofA_b}-\ref{sec:polynomial_precond}, a rigorous derivation of these methods, analyze their asymptotic stability and formal proximity to singularity. Then, in Section \ref{sec:numerical_experiments}, we present numerical experiments demonstrating their algorithmic efficiency on the 2D Laplacian and a system stemming from the Wilson-Dirac operator from lattice QCD; in the former, the shift-and-invert preconditioner allows for a reduction in the effective number of matrix-vector multiplications under large stiffness, while using Chebyshev polynomial preconditioning maintains them and substantially reduces the orthogonalization costs, and in the latter example Ritz-based polynomial preconditioning in combination with deflation leads to a substantial reduction in overall computational cost. The framework, specialized here to two types of preconditioners, can be in principle utilized for any function $f$, and it is particularly beneficial when one needs to evaluate $f(A)b_{i}, \ i = 1, 2, ..., n$, with $n$ relatively large.

\section{Related Work}
\label{sec:related_work}



The fundamental approximation of $f(A)b$ via projection onto low-dimensional Krylov subspaces was pioneered for symmetric matrices by Druskin and Knizhnerman \cite{druskin89} and extended to the general non-Hermitian Arnoldi framework by Saad \cite{saad92}. The theoretical foundation of these un-preconditioned methods was subsequently solidified by Hochbruck and Lubich \cite{hochbruck1997on}, who established sharp, rigorous convergence bounds for Krylov subspace approximations of the matrix exponential, demonstrating their superlinear convergence properties. A comprehensive treatment of matrix functions and their general evaluation can be found in the text by Higham \cite{higham2008functions}. 

While these standard approximations have been extensively studied \cite{saad, golub}, their convergence inevitably stagnates when $A$ exhibits an unfavorable spectral distribution or when $f(z)$ possesses singularities near the spectrum. Saad further investigated these limitations in his broader body of work on Krylov subspace behavior and polynomial approximations \cite{saad, saad92}, highlighting the critical need for acceleration techniques. Consequently, early pioneering work by Castillo and Saad \cite{saad97} explored preconditioning specifically for the matrix exponential operator, proposing both rational approximations and polynomial methods to enrich the Krylov subspace.

Loe and Morgan \cite{loe22} recently revisited polynomial preconditioning for GMRES in the context of linear systems, demonstrating that high-degree polynomials can significantly accelerate convergence for difficult matrices. Recent advancements in matrix-free approaches for $f(A)b$ have successfully utilized polynomial preconditioning and deflation techniques, particularly for the matrix inverse square root and the matrix sign function in lattice QCD \cite{ramirez23, frommer2024polynomial}. 

The specific interplay of dynamic Ritz-based polynomial preconditioning with exact LR-deflation formulated in our work is strongly inspired by the successful combination of these techniques for solving highly ill-conditioned linear systems on the coarsest grids of multigrid hierarchies in lattice QCD, as demonstrated by Espinoza-Valverde et al.\ \cite{espinoza2023coarsest}.

A critical bridge to adapting these techniques for general matrix functions is the concept of \textit{shift-invariance}. Expressing $f(A)b$ via the Cauchy contour integral naturally connects the evaluation of a matrix function to a continuum of shifted linear systems. The task of efficiently solving multi-shift linear systems using polynomial preconditioners was elegantly addressed by Ahmad, Szyld, and van Gijzen \cite{ahmad17}. 

Additionally, the efficiency of Krylov subspace methods for matrix functions is notoriously degraded when the matrix possesses eigenvalues close to the function's singularities. The resolution of this instability via augmented Krylov subspaces and LR-deflation was rigorously explored by Bloch et al.\ \cite{bloch07} for non-Hermitian overlap Dirac operators, which we take as the deflation approach to use in this work.

Our contribution unifies these domains. We lift the discrete shift-invariance framework of Ahmad et al.\ out of the realm of standard linear equations and inject it directly into the continuous contour integral definition of the Arnoldi method for $f(A)b$.

\section{The Arnoldi Method and Integral Representation}
\label{sect:Arnoldi_method_for_fofA_b}

We briefly re-derive the standard Arnoldi method for $f(A)b$ using integral representations \cite{saad92}. This serves as the foundation upon which our preconditioning framework is built.

For a function $f$ that is analytic on and inside a closed contour $\Gamma$ enclosing the spectrum of $A \in \mathbb{C}^{N \times N}$, one can express $f(A)b$ using the Cauchy integral formula:
\begin{equation} \label{eq:integral_form_fofAb}
    f(A)b = \frac{1}{2 \pi i}\int_{\Gamma} f(z) (zI_{N} - A)^{-1} b \, dz ,
\end{equation}
where $z$ is a complex variable, and $I_{N}$ is the identity matrix of dimension $N$.



The Arnoldi process generates an orthonormal basis $V_{m+1} = [v_1, v_2, \dots, v_{m+1}] \in \mathbb{C}^{N \times (m+1)}$ for the Krylov subspace $\mathcal{K}_m(A, b)$, satisfying the Arnoldi relation $AV_{m} = V_{m}H_{m} + h_{m+1,m}v_{m+1}e_{m}^{H}$.

By utilizing the Full Orthogonalization Method (FOM) and exploiting the \textit{shift-invariance} of Krylov subspaces---specifically the property that a shifted matrix generates the exact same Krylov space, $\mathcal{K}_{m}(zI_{N}-A, b) = \mathcal{K}_{m}(A, b)$ \cite{frommer1995cg}---we can approximate the solutions to the shifted linear systems inside the integral as $x_{m,z} = V_{m} (zI_{m}-H_{m})^{-1} V_{m}^{H} b$. Substituting this approximation back into the continuous integral gives the classic Arnoldi approximation for $f(A)b$:
\begin{equation}\label{eq:Arnoldi_method_for_fofA_b}
    f(A) b \approx \frac{1}{2 \pi i} \int_{\Gamma} f(z) V_{m} (zI_{m}-H_{m})^{-1} V_{m}^{H} b \, dz = V_{m} f(H_{m}) V_{m}^{H}b .
\end{equation}

\section{The General Preconditioned Framework}
\label{sec:proposed_framework}

We now extend the standard Arnoldi method to accommodate a preconditioner $\widehat{M}^{-1} \approx A^{-1}$. Our goal is to construct a single Krylov basis using the well-conditioned operator $A\widehat{M}^{-1}$ while preserving the integrity of the contour integral.

We begin by inserting a $z$-dependent preconditioner $\widehat{M}_{z}^{-1}$ into the shifted linear systems of the integral representation:
\begin{equation} \label{eq:integral_form_fofAb_preconditioned}
    f(A)b = \frac{1}{2 \pi i}\int_{\Gamma} f(z) \widehat{M}_{z}^{-1} \widehat{y}_{z} \, dz , \quad \text{where } (zI_{N} - A) \widehat{M}_{z}^{-1} \widehat{y}_{z} = b .
\end{equation}

To avoid generating a different Krylov subspace for every evaluation node $z \in \Gamma$, we rely on shift-invariance. Let $\widehat{M}_{0}^{-1}$ denote the base preconditioner. We desire the relation \cite{ahmad17}:
\begin{equation} \label{eq:krylov_subspace_shift_invariance}
    \mathcal{K}_{m}(A \widehat{M}_0^{-1}, b) = \mathcal{K}_{m}(\eta_{z}I_{N}-A\widehat{M}_0^{-1}, b) = \mathcal{K}_{m}((zI_{N}-A)\widehat{M}^{-1}_{z}, b) ,
\end{equation}
where $\eta_{z}$ is a scalar shifting parameter dependent on $z$. To satisfy this equality, it is sufficient to find an $\eta_z$ and $\widehat{M}_{z}^{-1}$ such that the operators are equivalent:
\begin{equation}\label{eq:relation_shift_inv_precond}
\eta_{z}I_{N} - A\widehat{M}_0^{-1} = (zI_{N} - A)\widehat{M}^{-1}_{z} .
\end{equation}

We generate the preconditioned Arnoldi relation using the base operator $A\widehat{M}_{0}^{-1}$:
\begin{equation} \label{eq:Arnoldi_preconditioned_M0}
 (A \widehat{M}_{0}^{-1}) V_{m} = V_{m} \widetilde{H}_{m} + h_{m+1,m}v_{m+1}e_{m}^{H} . 
\end{equation}

Applying the FOM condition to the equivalent preconditioned operator yields the subspace approximation $\widehat{y}_{z} \approx V_{m} (\eta_{z}I_{m} - \widetilde{H}_{m})^{-1} V_{m}^{H} b$. Substituting this into eq.\ \ref{eq:integral_form_fofAb_preconditioned} provides our fundamental preconditioned approximation, denoted as $f_m$:
\begin{equation}\label{eq:final_approximation_fm}
 f_{m} = \frac{1}{2 \pi i}\int_{\Gamma} f(z) \widehat{M}_{z}^{-1} V_{m} (\eta_{z}I_{m} - \widetilde{H}_{m})^{-1} V_{m}^{H} b \, dz .
\end{equation}

The subsequent sections explore how specific choices for $\widehat{M}_0^{-1}$ uniquely decouple eq.\ \ref{eq:final_approximation_fm}.

\section{Rational Shift-and-Invert Preconditioning}
\label{sec:shift_invert}

We first consider the case where the base preconditioner is a rational, shift-and-invert operator. 
\begin{assumption}
The matrix $A$ has a spectrum strictly contained in the right half of the complex plane, i.e., $\operatorname{Re}(\lambda) > 0$ for all $\lambda \in \sigma(A)$. This encompasses symmetric positive definite (SPD) matrices as well as non-Hermitian matrices with positive real parts.
\end{assumption}

Note that this assumption could be in principle relaxed, but we assume it to hold throughout this work. We define then the base preconditioner using a positive scalar shift $\mu > 0$:
\begin{equation}
    \widehat{M}_0^{-1} = (\mu I_N + A)^{-1}.
\end{equation}

\subsection{Closed-Form Approximation}
To satisfy the shift-invariance relation (eq.\ \ref{eq:relation_shift_inv_precond}), we equate the operators:
\begin{equation}
    \eta_{z}I_{N} - A(\mu I_N + A)^{-1} = (zI_{N} - A)\widehat{M}^{-1}_{z} .
\end{equation}
By defining the shifted preconditioner as $\widehat{M}_{z}^{-1} = \frac{\mu}{\mu + z}(\mu I_N + A)^{-1}$, simple algebraic manipulation isolates the required scalar shifting parameter:
\begin{equation}
    \eta_z = \frac{z}{\mu + z}.
\end{equation}

Substituting these explicit expressions into the fundamental approximation (eq.\ \ref{eq:final_approximation_fm}) yields the mapped integral formulation:
\begin{align}
 f_{m} &= \frac{1}{2 \pi i}\int_{\Gamma} f(z) \left( \frac{\mu}{\mu + z} \right) (\mu I_N + A)^{-1} V_{m} \left( \frac{z}{\mu + z} I_{m} - \widetilde{H}_{m} \right)^{-1} V_{m}^{H} b \, dz \\
 &= \frac{1}{2 \pi i} \mu (\mu I_N + A)^{-1} V_{m} \left[ \int_{\Gamma} f(z) \left( zI_{m} - (\mu + z)\widetilde{H}_{m} \right)^{-1} dz \right] (\beta e_1) . \label{eq:sai_integral_form}
\end{align}

To eliminate the numerical contour quadrature, we factor the term $(I_m - \widetilde{H}_m)$ out of the inverse inside the integral:
\begin{equation} \label{eq:resolvent_factored}
 \left( zI_{m} - (\mu + z)\widetilde{H}_{m} \right)^{-1} = (I_m - \widetilde{H}_m)^{-1} \left( zI_m - \mu(I_m - \widetilde{H}_m)^{-1}\widetilde{H}_m \right)^{-1} .
\end{equation}

By defining the transformed small subspace matrix $\widehat{H}_m := \mu(I_m - \widetilde{H}_m)^{-1}\widetilde{H}_m$, we can substitute eq.\ \ref{eq:resolvent_factored} back into the integral, allowing the $z$-independent terms to be extracted entirely. This yields an exact, closed-form extraction:
\begin{equation} \label{eq:closed_form_exact}
 f_m = \mu (\mu I_N + A)^{-1} V_m (I_m - \widetilde{H}_m)^{-1} f\left(\widehat{H}_m\right) (\beta e_1) .
\end{equation}

\subsection{Spectral Mapping and the Phantom Pole Phenomenon}
\label{sec:limitations}

The effectiveness of the rational preconditioner depends heavily on the choice of the shift parameter $\mu$. From the relation $\lambda_W = \frac{\lambda}{\mu + \lambda}$, we observe that $\mu$ dictates the spectral compression of the preconditioned operator $W = A(\mu I_N + A)^{-1}$. 

Consider the behavior of the required matrix inversion $(I_m - \widetilde{H}_m)^{-1}$ in the closed-form formulation (eq.\ \ref{eq:closed_form_exact}). As $\mu \to 0$, the largest eigenvalues $\lambda_W$ cluster near $1.0$. This inevitably forces the eigenvalues of the projected Hessenberg matrix $\widetilde{H}_m$ to approach $1.0$, rendering the matrix $(I_m - \widetilde{H}_m)$ highly singular. To shield the extraction step from this singularity, we establish a theoretical boundary. We require the distance between the largest eigenvalue of $W$ (corresponding to $\lambda_{\max}$ of $A$) and $1.0$ to be strictly greater than a tunable tolerance $\tau > 0$:
\begin{equation}
    1 - \lambda_W(\lambda_{\max}) = 1 - \frac{\lambda_{\max}}{\mu + \lambda_{\max}} = \frac{\mu}{\mu + \lambda_{\max}} > \tau.
\end{equation}
Assuming $\mu \ll \lambda_{\max}$, this provides a required lower bound for the shift:
\begin{equation}
    \mu > \lambda_{\max} \cdot \tau.
\end{equation}
The tolerance $\tau$ serves as a critical stabilization parameter, which can be dynamically varied to balance the preconditioning compression against the specific internal stiffness of the target function (e.g., the diffusion time step $c$ for the matrix exponential, see Section \ref{sec:laplacian_exp}). As $\tau$ is tuned down to accommodate increasingly difficult $f(A)b$ evaluations, the correspondingly smaller shift $\mu$ inherently reduces the diagonal dominance of the shifted linear systems $(\mu I_N + A)x = y$. While this renders the inner solves progressively more difficult for standard iterative methods, a properly constructed geometric or algebraic multigrid solver remains insensitive to these changes in $\mu$, retaining optimal convergence rates even as $\mu \to 0$.

\begin{figure}[htbp]
  \centering
  \includegraphics[width=0.65\textwidth]{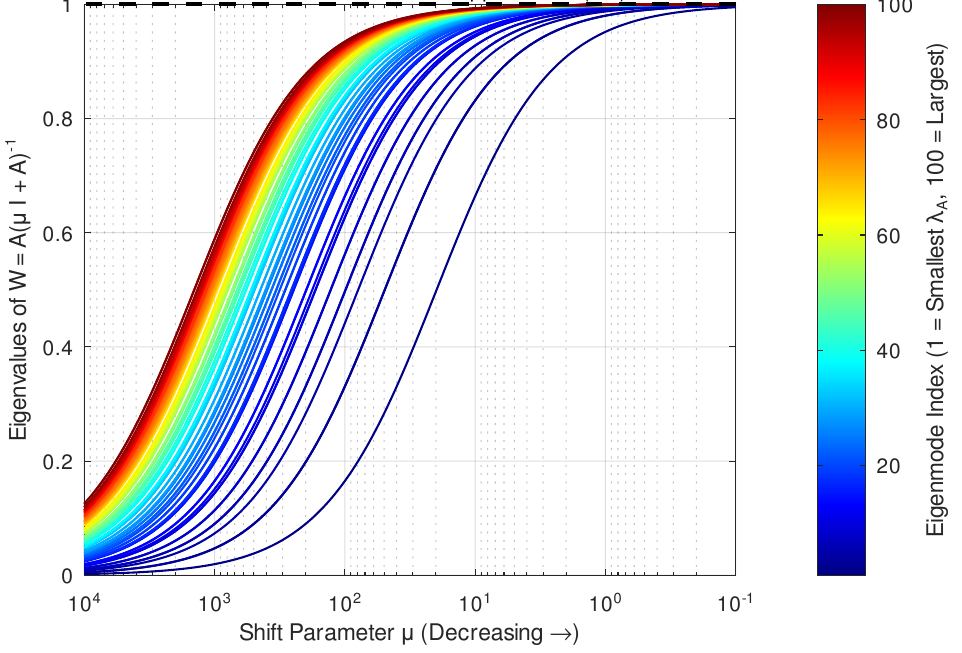}
  \caption{The 100 smallest eigenvalues of $W$, with $A$ the 2D Laplace matrix for a $512{\times}512$ grid, for which $\lambda_{min}=19.74$ and $\lambda_{max} = 2.1\cdot10^{6}$. When $\mu$ is forced to be relatively large (left side along the horizontal axis), the low modes are severely clustered towards zero. Shrinking $\mu$ (toward the right) cures this compression but forces the highest modes of $\widetilde{H}_{m}$ toward $1.0$.}
  \label{fig:varying_mu}
\end{figure}

\begin{figure}[htbp]
  \centering
  \includegraphics[width=0.65\textwidth]{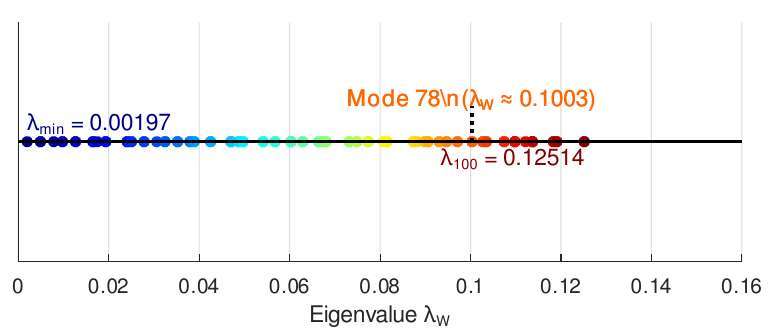}
  \caption{Distribution of the 100 smallest eigenvalues of $W$ for a large shift ($\mu=10^4$), for the same system as in Figure \ref{fig:varying_mu}. The extreme clustering near zero severely degrades the representational power of the Krylov subspace.}
  \label{fig:mu_distribution}
\end{figure}

However, fulfilling this theoretical condition introduces a severe secondary limitation. As visually evidenced in Figures \ref{fig:varying_mu} and \ref{fig:mu_distribution}, choosing a conservatively large $\mu$ (driven by a larger $\tau$) can considerably clusterize the low-frequency modes of $W$ toward zero. These near-null modes can render $\widehat{H}_{m}$ highly ill-conditioned, going against the main purpose of the preconditioning introduced here.


Therefore, we want to choose a relatively small $\mu$, associated with a small $\tau$.
This compresses the entire spectrum of $W$ into a tightly bounded cluster approaching $1.0$, which is theoretically ideal for rapid Krylov subspace convergence, while at the same time maintaining a relatively good conditioning of the matrices $\widetilde{H}_{m}$ and $\widehat{H}_{m}$.


This extreme clustering, though, introduces a finite-precision vulnerability. Because $W \approx I$, successive Krylov vectors generated by the Arnoldi process become nearly collinear. 
When orthogonalizing this basis using standard Modified Gram-Schmidt (MGS), the algorithm suffers from severe cancellation and loss of orthogonality. This introduces floating-point noise into the basis, producing spurious ``phantom poles" \cite{gonnet2013robust} (artificial Ritz values) near zero in the resultant matrix $\tilde{H}_m$.

To cure this instability, we enforce Double Modified Gram-Schmidt (DMGS) reorthogonalization during the preconditioned Arnoldi process. By explicitly executing a second projection pass at each iteration, we strip away the finite-precision noise before it can contaminate the subspace basis. Consequently, the phantom poles are largely suppressed, $\widetilde{H}_m$ stably and accurately represents the true, well-conditioned spectrum of $W$, and the exact closed-form extraction (eq. \ref{eq:closed_form_exact}) can be evaluated without requiring any external eigenvector deflation.

\section{Polynomial Preconditioning via Contour Quadrature}
\label{sec:polynomial_precond}

In case there is no multigrid implementation available for solving the large linear system $(\mu I_{N} + A)x=y$ required by rational methods, or if that approach is simply overkill, we introduce a polynomial preconditioner, $\widehat{M}_{0}^{-1} = p(A)$.

Before projecting onto a Krylov subspace, it is instructive to observe the exact continuous behavior of the preconditioned system. Recalling the exact integral representation and the shift-invariance relation (Eq. \ref{eq:relation_shift_inv_precond}) from Section \ref{sec:proposed_framework}, we can rewrite the shifted linear systems in a more compact form by introducing the shifted preconditioned operator:
\begin{equation} \label{eq:shifted_precond_op_def}
\widehat{A}_{z} := (zI_{N} - A) \widehat{M}_{z}^{-1} = \eta_{z}I_{N} - A\widehat{M}_{0}^{-1}.
\end{equation}
If we were to solve these preconditioned shifted linear systems exactly—meaning $\hat{y}_{z} = \widehat{A}_{z}^{-1} b$—we could inject this directly back into the exact contour integral, retrieving $f(A)b$ exactly:
\begin{equation}
f(A)b = \frac{1}{2\pi i}\int_{\Gamma} f(z) \widehat{M}_{z}^{-1} \widehat{A}_{z}^{-1} b \, dz.
\end{equation}

\subsection{Basis Selection: Monomial vs. Newton Formulation}
\label{sec:basis_selection}
To efficiently evaluate this framework, we must select an algebraic basis to represent the preconditioner $p(A)$. The choice of basis is not merely a formality; it dictates the fundamental numerical stability of both the continuous contour integral and the global matrix-vector assembly. We explicitly present two distinct formulations: the \textit{monomial basis}, reserved for low-degree simplicity, and the \textit{Newton basis}, which is mandatory for the numerical stability of high-degree applications.

\subsubsection{The Monomial Basis Formulation}
\label{sec:monomial_basis}
For benign operators requiring only minimal spectral compression, the preconditioner can be expressed as a low-degree polynomial $d$ in the standard monomial basis with known coefficients $\gamma_{j}$:
\begin{equation}
\widehat{M}_{0}^{-1} = \sum_{j=0}^{d}\gamma_{j}A^{j}.
\end{equation}
Following the shift-invariance relations established in Section \ref{sec:proposed_framework}, the shifted preconditioner similarly takes the form $\widehat{M}_{z}^{-1} = \sum_{j=0}^{d}\gamma_{z,j}A^{j}$, where the shifted coefficients are $\gamma_{z,j} = \sum_{\ell=j}^{d}\gamma_{\ell}z^{\ell-j}$ \cite{ahmad17}. Injecting this polynomial representation of $\widehat{M}_{z}^{-1}$ directly into our exact integral yields a summation of exact integrals, where the $z$-independent matrix powers $A^j$ can be safely factored out of the continuous quadrature:
\begin{equation}
f(A)b = \sum_{j=0}^{d} A^{j} \left[ \frac{1}{2\pi i}\int_{\Gamma}f(z) \left( \sum_{\ell=j}^{d}\gamma_{\ell}z^{\ell-j} \right) \widehat{A}_{z}^{-1} b \, dz \right]. \label{eq:exact_poly_integral_mono}
\end{equation}
By matching the equivalent operator definitions from Eq. \ref{eq:shifted_precond_op_def}, the continuous contour mapping scalar $\eta_z$ is identically evaluated via the coefficient sum:
\begin{equation} \label{eq:eta_mono}
\eta_{z} = \sum_{\ell=0}^{d}\gamma_{\ell}z^{\ell+1}.
\end{equation}

As demonstrated in our first numerical experiment (Section \ref{sec:laplacian_exp}), a low-degree ($d=5$) Chebyshev polynomial preconditions the 2D Laplacian. At such low degrees, the coefficients $\gamma_j$ remain small ($\mathcal{O}(10^1)$), and standard double-precision arithmetic easily absorbs the minimal coefficient growth, allowing for a straightforward evaluation via Horner's rule.

However, for challenging matrices requiring high polynomial degrees ($d \gtrsim 15$), the monomial coefficients of bounding polynomials (such as Chebyshev or Ritz-based) grow exponentially ($2^{d-1}$) with alternating signs. This triggers two severe numerical instabilities:
1. \textit{Contour Distortion:} Evaluating the scalar $\eta_z$ (Eq. \ref{eq:eta_mono}) via the monomial sum results in large floating-point cancellation, destroying the integration path and generating \texttt{NaN}s.
2. \textit{Intermediate Swell:} During the global full-space assembly of $f(A)b$ (Eq. \ref{eq:exact_poly_integral_mono}), Horner's rule scales intermediate vectors by $\mathcal{O}(10^{20})$ before canceling them down to an $\mathcal{O}(1)$ final vector, leading to a total loss of significance.

Consequently, while exact, the monomial basis is unviable for the high degrees necessitated by highly ill-conditioned cases.

\subsubsection{The Newton Basis Formulation}
\label{sec:newton_basis}
To cure the intermediate swell and catastrophic cancellation for highly ill-conditioned operators (e.g., the LQCD operator in Section \ref{sec:twisted_mass_experiment}, where $d \in [20, 100]$), we must abandon the monomial expansion. Instead, we formulate the framework directly in the factorized root representation using the Newton basis.

Given $d$ polynomial roots $\theta_i$ (either extracted dynamically as Ritz values or defined analytically as Chebyshev nodes), we define a sequence of stable basis operators:
\begin{equation}
    N_0(A) = I, \qquad N_j(A) = N_{j-1}(A)\left( I - \frac{A}{\theta_j} \right) = \prod_{i=1}^j \left( I - \frac{A}{\theta_i} \right).
\end{equation}
Applying this basis sequentially to a vector natively filters out the specific root $\theta_j$ at each step, bounding the vector norms and mitigating intermediate swell.

The shifted preconditioner $\widehat{M}_z^{-1}$ can be expressed as a linear combination of these stable operators:
\begin{equation}
    \widehat{M}_z^{-1} = \sum_{j=0}^{d-1} \beta_j(z) N_j(A), \qquad \text{where } \beta_j(z) = \frac{1}{\theta_{j+1}} \prod_{k=j+2}^d \left( 1 - \frac{z}{\theta_k} \right).
\end{equation}
Injecting this Newton basis into the exact integral, the $z$-independent stable operators $N_j(A)$ are similarly factored out:
\begin{equation}
f(A)b = \sum_{j=0}^{d-1} N_j(A) \left[ \frac{1}{2\pi i}\int_{\Gamma}f(z) \beta_j(z) \widehat{A}_{z}^{-1} b \, dz \right]. \label{eq:exact_poly_integral_newton}
\end{equation}
The contour shift scalar $\eta_z = z p(z)$ is now identically and stably evaluated via the factorized roots, bypassing any coefficient expansion:
\begin{equation}
    \eta_z = 1 - \prod_{k=1}^d \left( 1 - \frac{z}{\theta_k} \right).
\end{equation}
This root-based formulation protects both the continuous contour integration and the global matrix-vector assembly against finite-precision breakdown, enabling the stable deployment of arbitrarily high polynomial degrees.

Notably, deploying this Newton basis does not incur any additional computational overhead compared to the standard monomial expansion. Just as the monomial basis can be nested using Horner's rule to avoid computing explicit matrix powers, the Newton basis can be nested using a generalized Horner's rule \cite{higham2008functions}. If we want to compute the final full-space assembly $f_m = \sum_{j=0}^{d-1} N_j(A) w_j$ (where $w_j$ are the recombined decoupled vectors), we do not compute the operator $N_j(A)$ explicitly. Instead, because each basis term is defined recursively as $N_j(A) = N_{j-1}(A)\left( I - \frac{A}{\theta_j} \right)$, we can nest the factored roots from the inside out:
\begin{equation}
    f_m = w_0 + \left(I - \frac{A}{\theta_1}\right) \Bigg( w_1 + \left(I - \frac{A}{\theta_2}\right) \bigg( w_2 + \dots + \left(I - \frac{A}{\theta_{d-1}}\right) w_{d-1} \bigg) \Bigg).
\end{equation}
Algorithmically, this is evaluated by initiating a vector $v_{d-1} = w_{d-1}$, and then running a simple loop descending from $j = d-1$ down to $1$:
\begin{equation}
    v_{j-1} = w_{j-1} + \left(I - \frac{A}{\theta_j}\right) v_j.
\end{equation}
By distributing the vector $v_j$, the recursive update becomes:
\begin{equation}
    v_{j-1} = w_{j-1} + v_j - \frac{1}{\theta_j} A v_j.
\end{equation}
Looking closely at this final expression, there is exactly one application of the sparse operator $A$ per iteration. Because this loop runs exactly $d-1$ times, it costs exactly $d-1$ sparse matrix-vector multiplications (SpMVs). This matches the computational footprint of the monomial basis, achieving numerical stability without paying a single extra SpMV.

\subsection{The Illusion of the Closed Form}
Given that the shift-and-invert framework supports closed-form extraction, it is tempting to seek a closed-form solution for the polynomial base integrals using Cauchy's Residue Theorem. Theoretically, the integral evaluates to the sum of the residues at the roots of the equation $\eta_z - \lambda_r = 0$, where $\lambda_r \in \sigma(\tilde{H}_m)$.

However, this approach is numerically unstable. Root-finding in this context is highly sensitive to Wilkinson's phenomenon, where small floating-point errors cause the roots to scatter unpredictably across the complex plane.
Therefore, robust numerical quadrature along $\Gamma$ is required. Numerical quadrature avoids root finding and isolates interior branch cuts via the contour path.

\subsection{Krylov Projection and Schur Decomposition}
\label{sec:diagonalization}
Because evaluating the exact inverse $\widehat{A}_{z}^{-1}b$ at every complex quadrature node is computationally intractable, we must project the exact mathematical relations (Eqs. \ref{eq:exact_poly_integral_mono} and \ref{eq:exact_poly_integral_newton}) onto an Arnoldi Krylov subspace. As derived previously, the equivalent projected operator yields the subspace approximation $\widehat{A}_{z}^{-1}b \approx V_{m}(\eta_{z}I_{m} - \tilde{H}_{m})^{-1}V_{m}^{H}b$.

To establish a unified computational framework, let $B_j(A)$ denote the selected basis operator ($A^j$ or $N_j(A)$) and $\phi_j(z)$ denote its corresponding quadrature weight ($\sum \gamma_\ell z^{\ell-j}$ or $\beta_j(z)$). The fundamental preconditioned Krylov approximation becomes:
\begin{equation}
f_{m} = \sum_{j} B_j(A) V_{m} \left( \left[ \frac{1}{2\pi i}\int_{\Gamma}f(z) \phi_j(z) (\eta_{z}I_{m} - \tilde{H}_{m})^{-1} dz \right] \right) V_{m}^{H}b.
\end{equation}


By defining a set of small base integral matrices $G_{j}\in\mathbb{C}^{m\times m}$, performing numerical quadrature naively requires solving a dense linear system $(\eta_{z}I_{m} - \tilde{H}_{m})$ at every single complex node $z$ along $\Gamma$. To eliminate this full matrix inversion, we employ the stable \textit{Schur decomposition} \cite{davies2003schur}: $\tilde{H}_{m} = Q T Q^{H}$, where $T$ is upper triangular and $Q$ is unitary ($\kappa(Q)=1$).

Substituting this factorization allows $Q$ to be factored outside the integral. Because the final Krylov approximation only requires the action of these matrices on the canonical vector $e_{1}$, we define $c = Q^{H}e_{1}$. The required vector $g_{j} = G_j e_1$ becomes:
\begin{equation}
g_{j} = Q \underbrace{ \left[ \frac{1}{2\pi i}\int_{\Gamma}f(z) \phi_j(z) (\eta_{z}I_{m} - T)^{-1} c \, dz \right] }_{u_{j}} = Q u_{j}.
\end{equation}
At each quadrature node $z$, evaluating the integrand only requires solving the upper triangular linear system $(\eta_{z}I_{m} - T)v_{z} = c$ via standard back-substitution, reducing the dense matrix solve to a $\mathcal{O}(m^{2})$ operation.

\subsection{Avoiding Floating-Point Overflow: Dual-Regime Evaluation}
\label{sec:overflow_avoidance}
\ \\
\noindent When matrices arise from fine PDE discretizations, the evaluation nodes $|z|$ on $\Gamma$ can extend toward infinity. Evaluating the integrands directly will trigger unbounded floating-point growth, resulting in indeterminate NaNs (e.g., Inf / Inf) during the linear solve $(\eta_z I_m - T)^{-1} c$.

To prevent this transient overflow, we analytically factor $\eta_z$ out of the linear system, defining the inverted scalar $\xi_z := 1/\eta_z$. By factoring out $\eta_z$, we rewrite the shifted inverse applied to $c$ as $\eta_z^{-1} (I_m - \xi_z T)^{-1} c$. We merge this $\eta_z^{-1}$ factor with our basis weights $\phi_j(z)$ to define the scaled weights $\widetilde{\phi}_j(z) = \phi_j(z)/\eta_z$.

This rescaling provides unconditional stability at infinity, causing both $\xi_z$ and $\widetilde{\phi}_j(z)$ to safely underflow to $0$. However, this introduces a dual vulnerability when $|z|$ is very small (if $\Gamma$ passes near the origin), where evaluating $\xi_z$ triggers division-by-zero. Fortunately, when $|\eta_z|$ is small, the original unscaled system is well-conditioned. We therefore implement a dynamic, dual-regime evaluation:
\begin{itemize}
    \item \textit{Regime 1 (Small $z$):} If $|\eta_z| \leq 1$, compute $v_z = (\eta_z I_m - T)^{-1} c$ and use the unscaled integrand $f(z) \phi_j(z) v_z$.
    \item \textit{Regime 2 (Large $z$):} If $|\eta_z| > 1$, compute $\tilde{v}_z = (I_m - \xi_z T)^{-1} c$ and use the scaled integrand $f(z) \widetilde{\phi}_j(z) \tilde{v}_z$.
\end{itemize}

The specific formulation of the scaled weight $\widetilde{\phi}_j(z)$ depends on the chosen basis:
\begin{itemize}
    \item \textit{Monomial Scaling:} $\widetilde{\phi}_j(z) = \left( \sum_{\ell=0}^d \gamma_\ell z^{\ell+1-k} \right)^{-1}$ (where $k = \ell - j$). 
    \item \textit{Newton Scaling:} $\widetilde{\beta}_j(z) = \frac{\beta_j(z)}{\eta_z} = \frac{1}{\theta_{j+1}} \frac{ \prod_{k=j+2}^d (1 - z/\theta_k) }{ 1 - \prod_{k=1}^d (1 - z/\theta_k) }$. As $z \to \infty$, the factorized polynomial evaluation ensures that $\widetilde{\beta}_j(z)$ gracefully and exactly decays to zero, explicitly preventing overflow without requiring any monomial expansion.
\end{itemize}

\subsection{Formal Proximity to Singularity of the Resolvents}
The computational robustness of the polynomial framework fundamentally hinges on the distance to singularity of the inverted matrices in both regimes. For any matrix $M$, proximity to a singularity is governed by $1/||M^{-1}||$.

In Regime 1, the resolvent is $(\eta_{z}I_{m}-T)^{-1}$. Its proximity to singularity is tightly associated with the distance $\min_{w\in\sigma(T)}|\eta_{z}-w|$. When the evaluation node $z$ is small (e.g., passing near the origin), the shift parameter vanishes, $\eta_{z}\rightarrow0$. If the Hessenberg matrix $T$ possessed small Ritz values near the origin, this distance would approach zero, precipitating catastrophic numerical breakdown. However, as derived in Section \ref{sec:lr_deflation}, the integration of explicit LR-deflation robustly extracts critical low modes from the operator. Consequently, the spectrum $\sigma(T)$ can be physically depleted of small eigenvalues, formally guaranteeing that $|\eta_{z}-w|$ remains strictly bounded away from zero.

In Regime 2, the scaled resolvent is $(I_{m}-\xi_{z}T)^{-1}$. As the quadrature node $z$ extends along the contour away from $|z|=1$, the polynomial shift grows unboundedly $(|\eta_{z}|\rightarrow\infty)$ which conversely drives the factored scalar to zero: $\xi_{z}=1/\eta_{z}\rightarrow0$. As $\xi_{z}\rightarrow0$, the scaled matrix asymptotes to the identity matrix:
\begin{equation}
(I_{m}-\xi_{z}T)\Rightarrow I_{m}
\end{equation}
This proves that as the integration marches toward infinity, the distance to singularity for the scaled system formally approaches exactly 1, rendering the linear solve stable against overflow.

\subsection{LR-Deflation for the Polynomial Framework}
\label{sec:lr_deflation}
The efficiency of Krylov subspace approximations is severely degraded when the matrix possesses eigenvalues close to the integration contour \cite{bloch07}. In our polynomial framework, this manifests if the well-conditioned operator $A\hat{M}_{0}^{-1}$ still possesses critical low modes, causing the Ritz values within $\tilde{H}_{m}$ to land dangerously close to the poles of $(\eta_{z}I_{m}-\tilde{H}_{m})^{-1}$. We cure this instability by adapting the LR-deflation scheme directly to the preconditioned operator.

\subsubsection{The Spectral Splitting}
We specifically extract $p$ critical low eigenmodes of the preconditioned operator, $A\widehat{M}_0^{-1}$, rather than the original matrix $A$. Deflating the preconditioned operator directly accelerates the shifted linear systems being solved within the integral formulation. Furthermore, it is important to note that when the eigenvalues of $A$ are very small (near the origin), the polynomial preconditioner inherently behaves as a linear transformation (an order-1 polynomial). Consequently, the critical low modes of $A$ map closely to the low modes of $A\widehat{M}_0^{-1}$, making them an ideal target for extraction.

This extraction yields the right/left eigenvectors $R_{p}$, $L_{p}\in \mathbb{C}^{N\times p}$ such that $(A\hat{M}_{0}^{-1})R_{p}=R_{p}\Lambda_{p}$ and $L_{p}^{H}(A\hat{M}_{0}^{-1})=\Lambda_{p}L_{p}^{H}$, normalized as $L_{p}^{H}R_{p}=I_{p}$. We decompose the initial vector $b$ into a critical component $b_{\parallel}=R_{p}(L_{p}^{H}b)$ and a deflated remainder $b_{\ominus}=b-b_{\parallel}$.

\subsubsection{The Exact Deflated Component}
Because $b_{\parallel}$ resides within the invariant subspace spanned by $R_{p}$, the application of the matrix function on this component can be evaluated exactly \cite{bloch07}, completely bypassing the continuous contour quadrature.

Since $A$ and the preconditioned operator $A\widehat{M}_0^{-1}$ commute, they fundamentally share a common eigenbasis. However, because the polynomial preconditioner aggressively folds and clusters the spectrum, it introduces severe eigenspace degeneracy in $A\widehat{M}_0^{-1}$. The basis vectors $R_p$ returned by standard eigensolvers (e.g., ARPACK) merely span this degenerate plane and do not necessarily mutually diagonalize the original matrix $A$. Consequently, extracting only the diagonal elements is insufficient, as it discards critical off-diagonal cross-terms that represent finite-precision basis rotations and near-degenerate mode mixing.

To evaluate the exact deflated solution, we must form the exact dense projection of $A$ onto this subspace. For the general non-Hermitian case, we utilize the bi-orthogonal left and right eigenvectors ($L_p^H R_p = I_p$) to perform an oblique projection:
\begin{equation}
    H_p = L_p^H A R_p .
\end{equation}
In the Hermitian case, the eigenvectors naturally form an orthogonal basis ($V_p^H V_p = I_p$), and this simplifies to the standard orthogonal projection $H_p = V_p^H A V_p$.

The exact action of the matrix function on the deflated space evaluates directly on this small dense matrix $H_p$, automatically capturing the true physics of the operator regardless of degeneracy. For the matrix inverse square root $f(A) = A^{-1/2}$, this corresponds to finding $y = H_p^{-1/2} \rho$, where the mapped right-hand side is $\rho = L_p^H b$.

To maximize numerical stability when evaluating this dense matrix function, we avoid explicit eigendecompositions or unstable matrix inversions (\texttt{inv()}). Instead, we evaluate the principal matrix square root $S_p = \mathrm{sqrtm}(H_p)$ using the stable Schur method \cite{bjorck1983schur, higham1987computing}. We subsequently solve the linear system $S_p y = \rho$ directly via backward/forward substitution (utilizing the standard backslash \texttt{\textbackslash} operator). The exact full-space deflated solution is then mapped back and assembled as:
\begin{equation}
    f_{\mathrm{exact\_def}} = R_p y .
\end{equation}

Because it is evaluated analytically in closed form, this exact spatial component is decoupled from the iterative Krylov procedure and is simply added as a constant vector to the final full-space assembly.

\subsubsection{The Deflated Krylov Component and Recombination}
For the remainder vector $b_{\ominus}$, we execute the Arnoldi process on $A\hat{M}_{0}^{-1}$ starting with $v_{1}=b_{\ominus}/||b_{\ominus}||_{2}$. Because $L_{p}^{H}b_{\ominus}=0$, the generated Krylov subspace is prevented from mixing with the critical eigendirections. The resulting Hessenberg matrix $\tilde{H}_{m}$ is naturally depleted of unstable low modes. We recombine these decoupled branches in the final algorithmic implementation.




\subsection{A Matrix-Free Stopping Criterion}\label{sec:stop_criterium_POLY_DEF}
In the polynomial framework, explicitly evaluating the full-space vector $f_{m}$ requires running $d$ computationally expensive sparse matrix-vector multiplications (SpMVs). Constructing $f_{m}$ at every Arnoldi step merely to check for convergence would severely degrade algorithmic efficiency.

To circumvent this, we track the stagnation of a single polynomial component directly in the small projected space. The base term corresponding to $j=0$ defines the identity operator in both bases ($A^0 = I$ and $N_0(A) = I$). Let $y_{m,0}=Q\zeta_{0}(T,\mathrm{kry})$ denote the small-space coefficient vector of length $m$ for this base degree term. 

When utilizing explicit LR-deflation alongside the polynomial preconditioner, we track convergence by monitoring the relative stagnation of exclusively the Krylov-projected component, $w_{m,0,\mathrm{kry}}=V_{m}y_{m,0}$. We track the relative change of this component:
\begin{equation}
\epsilon_{m}=\frac{||w_{m,0,\mathrm{kry}}-w_{m-1,0,\mathrm{kry}}||_{2}}{||w_{m,0,\mathrm{kry}}||_{2}}
\end{equation}
Because the Arnoldi basis $V_{m-1}$ constitutes the first $m-1$ columns of $V_{m}$, we can equivalently pad the previous coefficient vector with a zero, yielding $V_{m-1}y_{m-1,0}=V_{m}(y_{m-1,0}^T, 0)^T$. Substituting this into the relative change expression and utilizing the exact isometry of the orthonormal basis $V_{m}$ yields:
\begin{equation} \label{eq:stop_criter_POLY_DEF}
\epsilon_{m}=\frac{||V_{m}(y_{m,0}-(y_{m-1,0}^T, 0)^T)||_{2}}{||V_{m}y_{m,0}||_{2}}=\frac{||y_{m,0}-(y_{m-1,0}^T, 0)^T||_{2}}{||y_{m,0}||_{2}}< \mathrm{tol}.
\end{equation}
This derivation provides an exact, scale-invariant relative stopping criterion for the base polynomial component of the Krylov subspace. It can be evaluated at every iteration utilizing exclusively the small coefficient vectors, demanding absolutely zero SpMVs, while the explicitly assembled exact deflated solution is safely ignored during the iterative tracking phase.

\section{Numerical Experiments} 
\label{sec:numerical_experiments}
To empirically validate the efficiency, stability, and convergence of the proposed preconditioning frameworks, we evaluate their performance across two distinct benchmark problems. The primary goal is to demonstrate the algorithmic scalability and robustness when approximating $f(A)b$ under conditions of extreme PDE stiffness and unfavorable spectral distributions. Although all the methods outlined until now can be used for both Hermitian and non-Hermitian systems, the two numerical experiments in the next section deal with Hermitian a Hermitian matrix $A$.

All numerical experiments were performed serially on a 13-inch Apple MacBook Air equipped with an Apple M3 processor and 24 GB of RAM, running macOS 14.6. The algorithms were implemented and executed in GNU Octave version 11.1.0. The Octave code for generating the results presented here is available on GitHub\footnote{\url{https://github.com/Gustavroot/precFunDefl}}.

\subsection{Experimental Setup and Methodologies}
For both experiments, we introduce three algorithmic approaches:
\begin{enumerate}
    \item \textit{Traditional Arnoldi (STD-ARN):} The baseline Krylov subspace method without preconditioning.
    \item \textit{Shift-and-Invert Preconditioning (SAI-MG):} The algorithm derived in Section \ref{sec:shift_invert}, evaluating the closed-form formulation with reorthogonalization to preserve internal numerical stability against inner-solver noise. The internal sparse linear system solutions required for $(\mu I_{N}+A)^{-1}$ are executed using a multigrid solver.
    \item \textit{Contour-Decoupled Polynomial Method (POLY-DEF):} The algorithm developed in Section \ref{sec:polynomial_precond}, utilizing Schur factorization, LR-deflation on the preconditioned operator, and dual-regime contour quadrature.
\end{enumerate}

As detailed in the following sections, SAI-MG is extensively evaluated on the 2D Laplacian. Conversely, for the twisted-mass Wilson-Dirac operator evaluated in Section \ref{sec:twisted_mass_experiment}, while the $N=32768$ problem size and its extreme condition number would theoretically justify a multigrid approach, a tailored multigrid solver for this specific formulation is not available in our testing environment. Developing a dedicated multigrid solver solely for this matrix falls outside the scope of this work. Thus, we focus on the POLY-DEF method for this particular case.

\subsection{The Matrix Exponential and the 2D Laplacian}
\label{sec:laplacian_exp}
In the first experiment, we compute the action of the matrix exponential, $f(A)=\exp(-cA)$ for a positive scalar c, representing the formal solution to a parabolic diffusion PDE. We construct $A\in\mathbb{R}^{N\times N}$ as the standard 5-point finite difference discretization of the 2D Laplacian operator on a uniform Cartesian grid. For initial spectral visualizations, we utilize smaller grids, while the primary convergence evaluations are executed on a $512\times512$ grid $(N=262,144)$. We start with POLY-DEF in the next section to illustrate how dot products are exchanged for matrix-vector multiplications in that method, and then by moving to SAI-MG we can actually aspire to an overall reduction in both Krylov dimension and matrix-vector multiplications.

\subsubsection{Contour-Decoupled Polynomial Method (POLY-DEF)}
\begin{figure}[htbp]
    \centering
    \includegraphics[width=0.85\textwidth]{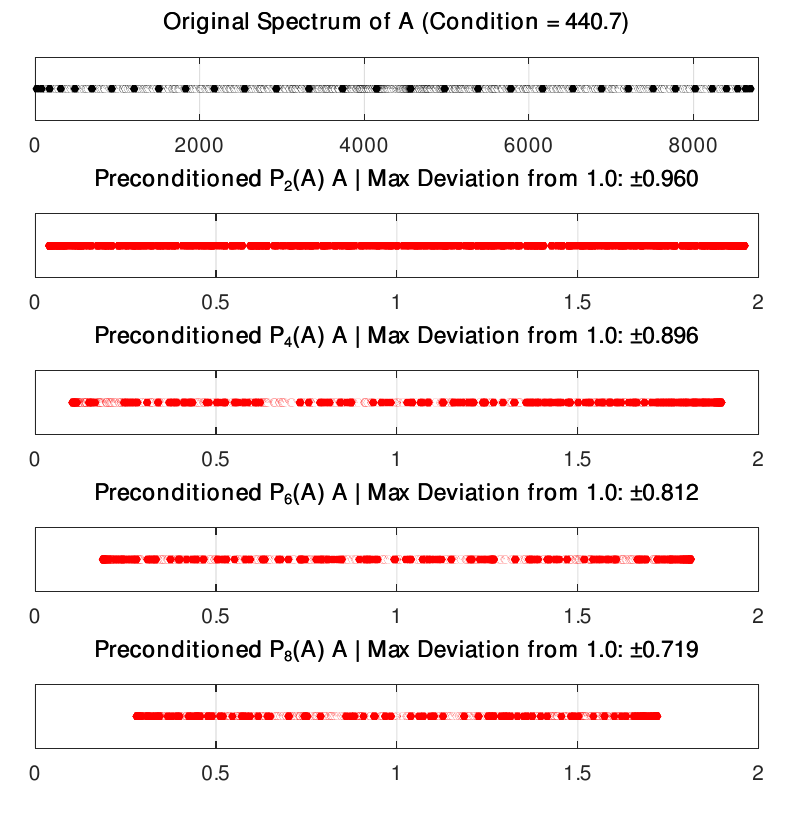}
    \caption{Original spectrum of A (top) and the compressed spectra of the preconditioned operator $P_{d}(A)A$ for varying polynomial degrees $d\in\{2,4,6,8\}$. The highly spread eigenvalues are aggressively clustered around 1.0, drastically reducing the condition number and mitigating PDE stiffness prior to the Krylov projection. The matrix $A$ corresponds to a 2D Laplacian with a $32\times32$ grid.}
    \label{fig:spectral_squeeze}
\end{figure}

\paragraph{Spectral Compression} Before evaluating the full Krylov convergence, we first visually demonstrate the spectral compression power of the Chebyshev polynomial preconditioner. Figure \ref{fig:spectral_squeeze} illustrates the effect of the preconditioned operator $M=P_{d}(A)A$ acting on a 2D Laplacian with a $32\times32$ grid $(N=1024)$. The original matrix exhibits a large condition number $(\kappa\approx440.7)$ with its spectrum spread extensively along the positive real axis.

By applying the matrix-free preconditioner, the minimax property of the Chebyshev polynomial compresses the scattered eigenvalues into a tightly bounded cluster centered at 1.0. As the polynomial degree $d$ is increased from 2 to $8,$ the maximum spectral deviation from 1.0 is reduced to $\pm0.719$. This extreme clustering transforms the initially ill-conditioned PDE operator into a well-conditioned one, forming the foundation for the rapid convergence of the decoupled Arnoldi framework.

\paragraph{Matrix-Free Error Tracking} To empirically validate the reliability of the matrix-free stopping criterion derived in eq.\ \ref{eq:stop_criter_POLY_DEF}, we tracked the internal small-space stagnation against the true full-space relative error, $||f_{m}-f_{\mathrm{exact}}||_{2}/||f_{\mathrm{exact}}||_{2}$. For this validation, we utilized a moderately ill-conditioned $64\times64$ Laplacian grid $(N=4096)$ with $c=100/\lambda_{\max}$ and computed the exact dense matrix exponential $f_{\mathrm{exact}}$ explicitly.

\begin{figure}[htbp]
    \centering
    \includegraphics[width=0.65\textwidth]{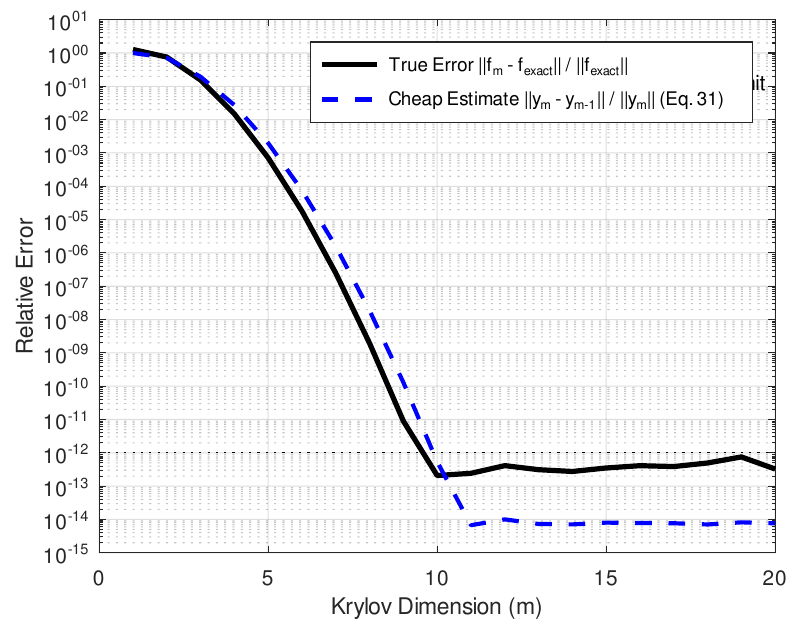}
    \caption{Comparison of the true full-space relative error (solid black line) and the cheap matrix-free stagnation estimate (dashed blue line) for the POLY-DEF algorithm $(d=5)$. The internal tracking models the descent rate of the true error, validating its use as a highly efficient stopping criterion. The matrix $A$ corresponds to a 2D Laplacian with a $32\times32$ grid and $c = 100/\lambda_{max}$.}
    \label{fig:error_proxy}
\end{figure}

As shown in Figure \ref{fig:error_proxy}, the proposed contour-decoupled polynomial framework drives the relative error down to the desired tolerance $(\approx10^{-12})$ in merely 10 Krylov iterations. The cheap relative stagnation estimate, evaluated exclusively on the base polynomial component in the small Schur space, shadows the true full-space error descent rate. This confirms that the algorithm can confidently halt execution dynamically without incurring the heavy SpMV penalty of Horner's full-space reconstruction at each iterative step.

\paragraph{Algorithm Convergence} Because $f(z)=\exp(-cz)$ is an entire function, it possesses no branch cuts. We select $\Gamma$ as a standard parabolic contour opening to the right half-plane \cite{weideman2007parabolic,schmelzer2007evaluating}. This encloses the positive real spectrum of the discretized Laplacian. The rapid asymptotic decay of the integrand eliminates strict requirements on tracking contour endpoints.

For POLY-DEF, the preconditioner degree $d$ is tuned empirically based on the spectral analysis above. Because the Laplacian lacks severe spectral outliers near the origin, deflation is less critical here. Consequently, a low degree of $d=5$ allows for the direct, efficient deployment of the monomial basis evaluated using standard double-precision Horner's rule, safely utilizing the monomial dual-regime scaling (Section \ref{sec:overflow_avoidance}) to robustly avoid transient overflow when mapping the highest Laplacian eigenvalues. We note that, as we are using Chebyshev polynomials in this example, we could have chosen to evaluate $\eta_z$ and the overall polynomial $f_m$ in a more numerically stable manner (either in factorized form or via Clenshaw's algorithm), but we maintain the monomial basis in this first example for algorithmic simplicity to demonstrate its viability at low degrees.

To evaluate the convergence under extreme stiffness, we scaled the problem to a $512\times512$ grid $(N=262,144)$ and set the diffusion time step to $c=10,000/\lambda_{\max}$. Because the discretized Laplacian is symmetric, the Arnoldi projection formally reduces to the short-recurrence Lanczos process. However, constructing the final full-space solution $f_{m}$ still traditionally requires storing the entire large Krylov basis $V_{m}$ in memory. To achieve a similarly low memory footprint as our matrix-free polynomial method, the baseline would need to employ a two-pass Lanczos approach, which entirely avoids storing the basis but doubles the required number of SpMVs. While we utilize the standard single-pass STD-ARN as our baseline to verify exact algebraic convergence and establish a conservative SpMV floor, it is critical to note that a memory-equivalent two-pass baseline would perform even worse in terms of raw operator applications.



Under these conditions, the standard single-pass STD-ARN struggles with the spectral spread. It requires 600 sparse matrix-vector multiplications (SpMVs) and an extensive Krylov subspace of dimension $m=600$ to reach the desired tolerance.

\begin{figure}[htbp]
    \centering
    \includegraphics[width=0.65\textwidth]{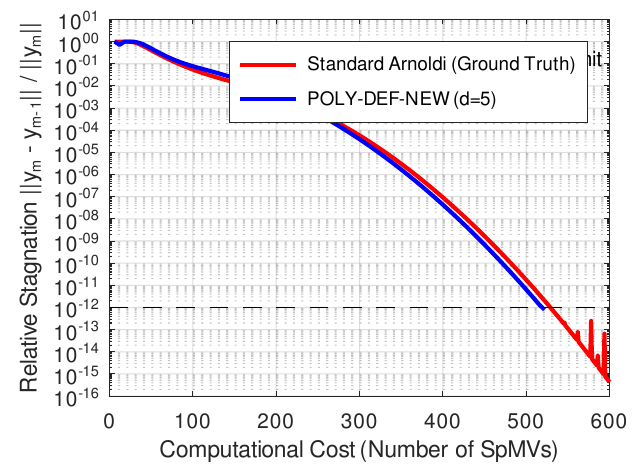}
    \caption{Convergence comparison on the highly stiff $512\times512$ Laplacian grid $(N=262,144$, $c=10000/\lambda_{\max})$. The proposed POLY-DEF method (blue) matches the precision of the standard Arnoldi method (red) using fewer SpMVs and drastically reduces the Krylov subspace dimension and its associated orthogonalization costs.}
    \label{fig:laplacian_convergence}
\end{figure}

In contrast, the POLY-DEF method (with polynomial degree $d=5$) overcomes the high-frequency stiffness. It converges to the precision floor in $\approx88$ Krylov iterations. Although each iteration requires $d+1$ SpMVs (totaling 527 SpMVs), the reduction in the Krylov subspace dimension eliminates the quadratic orthogonalization bottleneck. Figure \ref{fig:laplacian_convergence} illustrates this dynamic, plotting the internal relative stagnation against the true computational SpMV cost.

\subsubsection{Exact Shift-and-Invert Preconditioning (SAI-MG)}
Having established the efficiency of the matrix-free polynomial approach, we now turn our attention to the exact rational preconditioning (SAI-MG). For this method, the internal sparse linear systems are executed using a highly optimized geometric multigrid (GMG) solver. To safeguard convergence against inner-solver noise while isolating the evaluation from the unbounded contour, we deploy the exact closed-form formulation derived in Section \ref{sec:limitations} utilizing Double MGS reorthogonalization.

To evaluate the algorithmic scalability of SAI-MG, we examine the highly stiff $512\times512$ Laplacian $(N=262,144)$ evaluated under increasingly difficult diffusion time steps c.
To compensate for this increased difficulty within the SAI-MG framework, we dynamically tune the threshold parameter $\tau$ to decrease the preconditioner shift $\mu=\lambda_{\max}\cdot\tau$ This pushes the active spectrum of $W$ deeper toward 1.0, aggressively compressing the eigenvalues to systematically offset the increased stiffness of $f(A)$. Importantly, while shrinking $\tau$ technically reduces the diagonal dominance of the shifted operator $(\mu I_{N}+A)$, the GMG solver remains fundamentally insensitive to this variation, seamlessly preserving optimal $\mathcal{O}(N)$ convergence rates across all inner solves.

\begin{figure}[htbp]
    \centering
    \begin{subfigure}[b]{0.32\textwidth}
        \includegraphics[width=\linewidth]{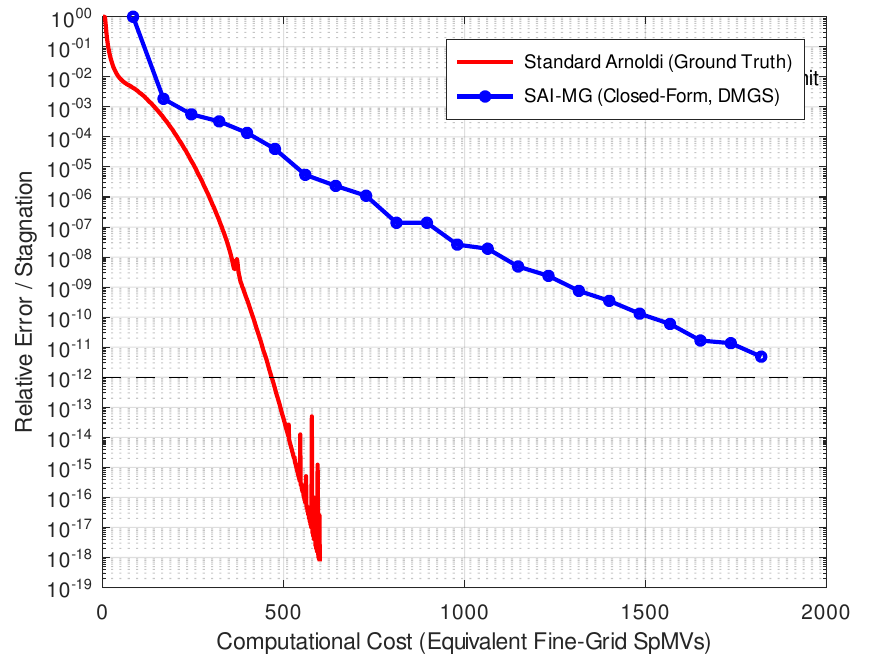}
        \caption{Case 1: $c_{1}=10^4$,\\$\mu=\lambda_{max}/10^3$}
    \end{subfigure}
    \hfill
    \begin{subfigure}[b]{0.32\textwidth}
        \includegraphics[width=\linewidth]{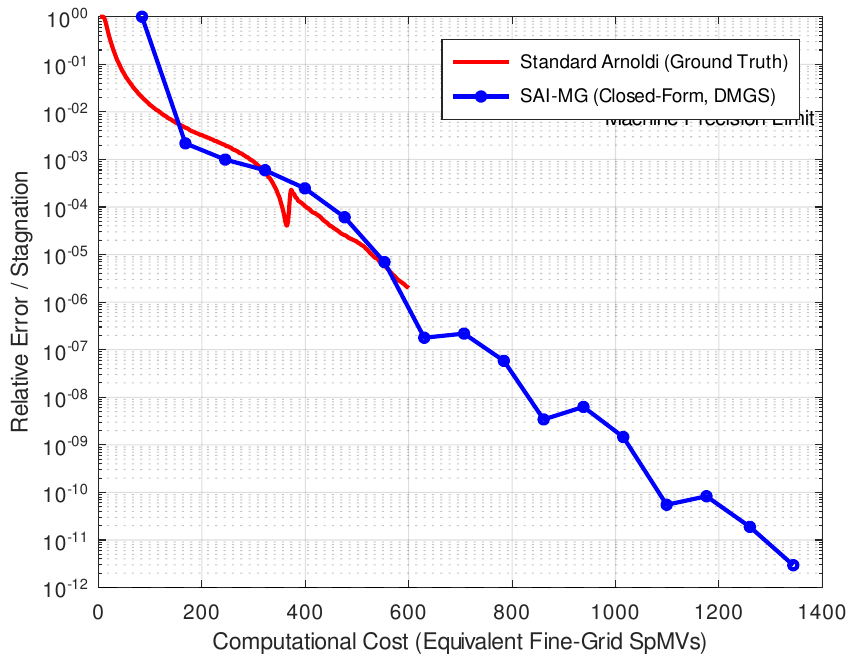}
        \caption{Case 2: $c_{2}=10^5$,\\$\mu=\lambda_{max}/10^4$}
    \end{subfigure}
    \hfill
    \begin{subfigure}[b]{0.32\textwidth}
        \includegraphics[width=\linewidth]{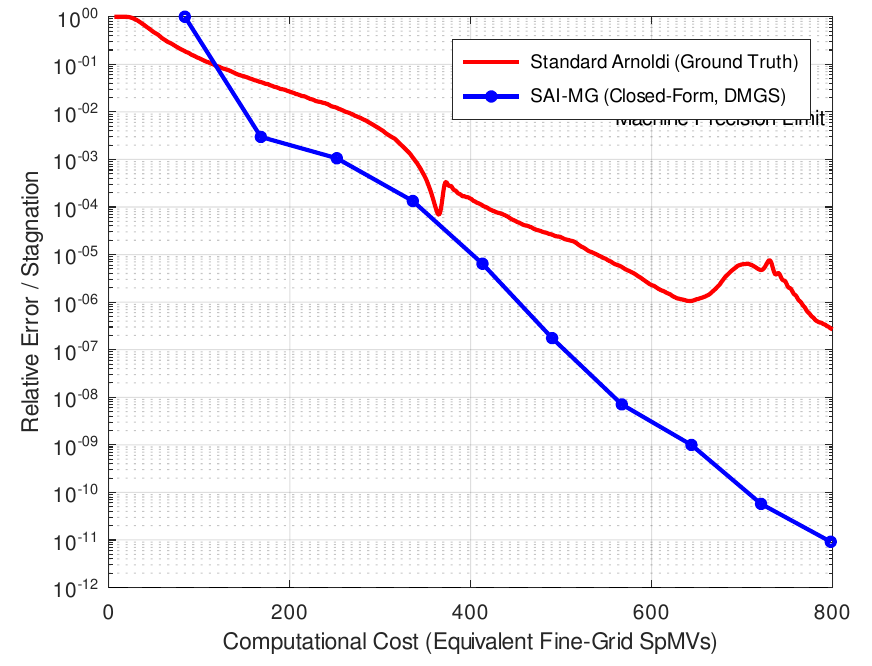}
        \caption{Case 3: $c_{3}=10^6$,\\$\mu=\lambda_{max}/10^5$}
    \end{subfigure}
    \caption{Convergence of SAI-MG vs. STD-ARN for increasingly stiff diffusion parameters $c$. As $c$ increases, the tolerance $\tau$ is dynamically tuned to shrink the shift $\mu$, allowing SAI-MG to absorb the stiffness and maintain rapid convergence while the baseline flatlines. The matrix $A$ corresponds to a 2D Laplacian with a $512\times512$ grid.}
    \label{fig:laplace_saimg_conv}
\end{figure}

Figure \ref{fig:laplace_saimg_conv} illustrates this convergence behavior across three distinct scales of extreme stiffness.
As the diffusion time step increases, the un-preconditioned STD-ARN baseline rapidly degrades, ultimately stagnating due to the large spectral gap.
Conversely, SAI-MG smoothly absorbs the stiffness through the $\tau$-tuned preconditioned operator, reaching the desired tolerance in a smaller number of equivalent fine-grid SpMVs. It is critical to emphasize that this efficiency gap is, in reality, significantly wider than plotted: to maintain a comparable $\mathcal{O}(1)$ -vector memory footprint to SAI-MG, the STD-ARN baseline would be practically forced to utilize a two-pass Lanczos algorithm, which doubles its accumulated SpMV costs.

\begin{figure}[htbp]
    \centering
    \begin{subfigure}[b]{0.32\textwidth}
        \includegraphics[width=\linewidth]{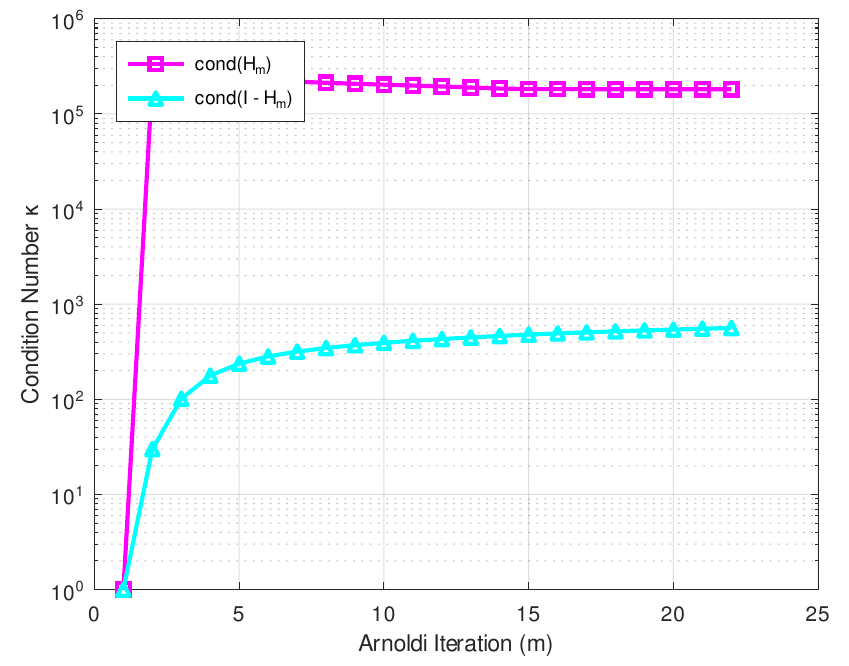}
        \caption{Case 1: $c_{1}$}
    \end{subfigure}
    \hfill
    \begin{subfigure}[b]{0.32\textwidth}
        \includegraphics[width=\linewidth]{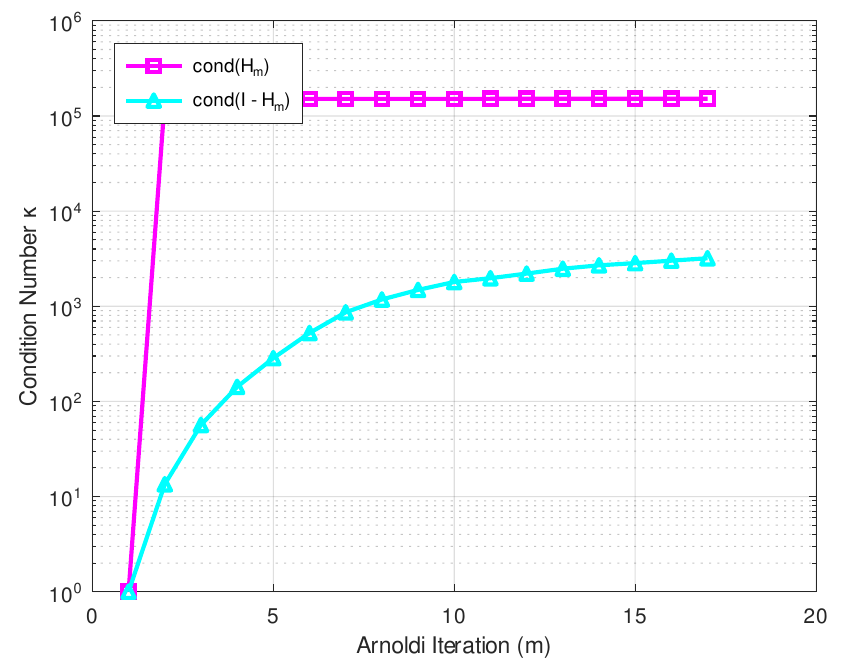}
        \caption{Case 2: $c_{2}$}
    \end{subfigure}
    \hfill
    \begin{subfigure}[b]{0.32\textwidth}
        \includegraphics[width=\linewidth]{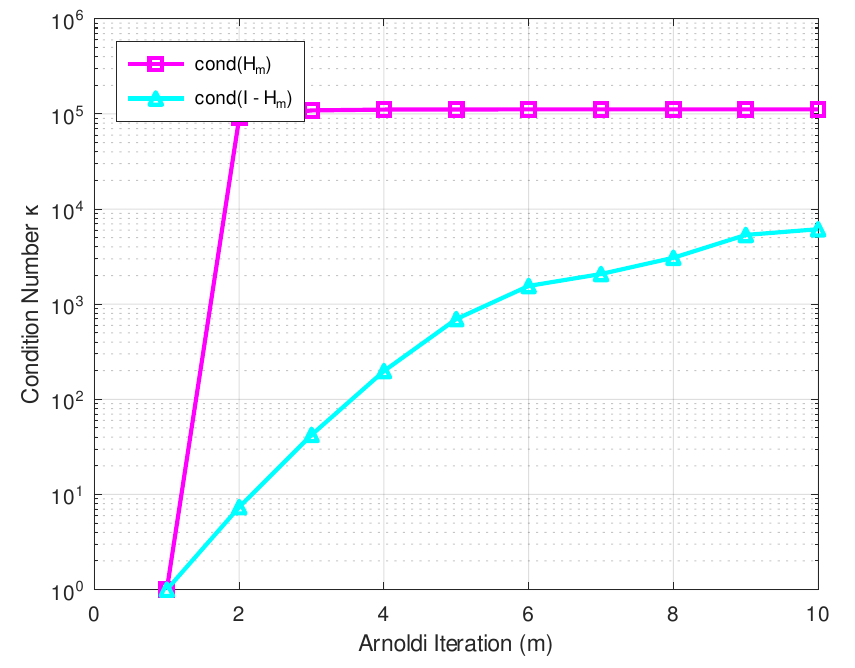}
        \caption{Case 3: $c_{3}$}
    \end{subfigure}
    \caption{Conditioning history of $\tilde{H}_{m}$ and $(I_{m}-\tilde{H}_{m})$ corresponding to the runs in Figure \ref{fig:laplace_saimg_conv}. The implementation of DMGS reorthogonalization successfully prevents the formation of phantom poles, maintaining structural stability inside the small-space extraction.}
    \label{fig:laplace_saimg_cond}
\end{figure}

Furthermore, Figure \ref{fig:laplace_saimg_cond} verifies the internal stability of the closed-form extraction sequence. Thanks to the Double MGS reorthogonalization implemented inside the preconditioned Arnoldi loop, the catastrophic formation of finite-precision phantom poles is entirely suppressed. The tracked condition numbers $\kappa(\tilde{H}_{m})$ and $\kappa(I_{m}-\tilde{H}_{m})$ remain safely bounded throughout the iteration sequence, verifying that the small-space dense matrix evaluation $f(\hat{H}_{m})$ can safely proceed without ever requiring an external implementation of exact eigenvector deflation.

\subsection{The Matrix Inverse Square Root on a Multigrid Twisted Mass Operator}
\label{sec:twisted_mass_experiment}
While the matrices evaluated in the previous sections validate the foundational mechanics of the proposed preconditioning frameworks, their spectral distributions are relatively benign. Low-degree polynomials can readily interpolate across their spectra, meaning they do not sufficiently stress-test the preconditioning framework in the presence of severe singularities. To robustly illustrate the efficacy of the POLY-DEF algorithm—and particularly the critical necessity of the Newton root basis for high degrees—we must evaluate an operator that exhibits the extreme spectral pathologies characteristic of realistic, large-scale scientific computations.

To this end, we extract a highly ill-conditioned, non-Hermitian operator arising from a state-of-the-art multigrid hierarchy used in Lattice Quantum Chromodynamics (QCD) \cite{gattringer2009quantum}. Specifically, we consider the maximally twisted mass fermion formulation, generated using the DD$\alpha$AMG solver \cite{frommer2013adaptive,simone2019simulating}. The original fine-grid ensemble corresponds to a large $96\times48^{3}$ lattice tuned to the physical pion mass. Because evaluating the exact inverse square root $f(A)b=A^{-1/2}b$ on this finest grid is computationally prohibitive for algorithmic benchmarking on standard workstations, we utilize the operator constructed at the coarsest level of the multigrid hierarchy, yielding a sparse non-Hermitian matrix of dimension $N=32768$.

The original coarse-grid twisted mass operator $D$ possesses a smallest eigenvalue of $1.4425\times10^{-2}$.
The corresponding Hermitian Positive Definite (HPD) base operator evaluated in our framework, $A=DD^{\dagger}$, is known to be bounded from below by $1.0{\cdot}10^{-4}$, possessing a minimum eigenvalue of $\lambda_{\min}\approx1.0680\times10^{-4}$ in our particular matrix. From a numerical linear algebra perspective, this coarsest-level matrix serves as an ideal benchmark: it inherits and mimics the extreme spectral difficulties of the finest-level operator, while operating at a computationally tractable dimension.

\begin{figure}[htbp]
    \centering
    \includegraphics[width=0.65\textwidth]{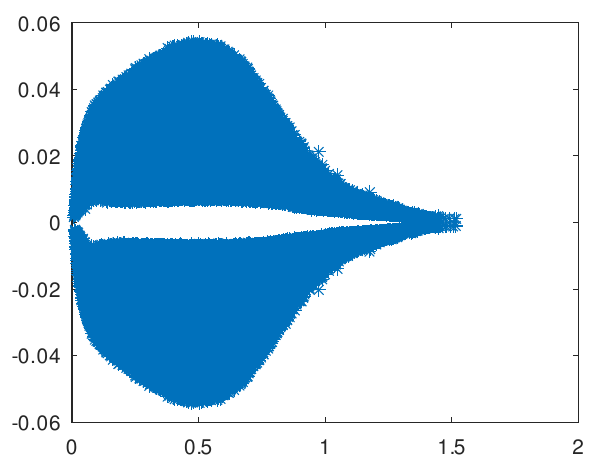}
    \caption{Eigenvalue spectrum of the coarsest-level non-Hermitian twisted mass operator ($N=32768$). The complex eigenvalues exhibit a severe, dense accumulation of low-frequency modes near the origin, presenting a critical convergence barrier for standard Krylov methods.}
    \label{fig:tm_coarse_spectrum}
\end{figure}

The defining pathological feature of this twisted mass operator is a uniquely high density of low-frequency modes accumulating near the origin (as illustrated in Figures \ref{fig:tm_coarse_spectrum} and \ref{fig:chebyshev_spectra}). This dense accumulation is a physical product of both the operator's physical reality and the twisted mass term itself. When evaluating the inverse square root $f(z)=z^{-1/2}$, this dense cluster of eigenvalues sits perilously close to the branch cut singularity at $z=0$, damaging convergence of Krylov-based methods in general.


By employing the POLY-DEF framework, we address this bad conditioning. We apply exact LR-deflation to extract the critical low-mode invariant subspace of the preconditioned operator, evaluating it via the decoupled scalar contour quadrature (as derived in Section \ref{sec:lr_deflation}). This depletes the dense near-origin cluster from the active spectrum of the operator. The contour-decoupled polynomial preconditioner is then deployed on the deflated remainder, compressing the well-separated, high-frequency spectrum towards unity without interference from the singularity. By decoupling the severe low-frequency modes from the bulk spectrum, the algorithm overcomes the low-mode stagnation and restores the rapid convergence characteristic of preconditioned Krylov methods.

\subsubsection{Deflated Polynomial Method (POLY-DEF)}
The operator at hand fundamentally tests the stability limits of Krylov methods due to the wide spread of complex eigenvalues and a high density of small eigenmodes near the origin. We expect STD-ARN to stall severely.

For the polynomial preconditioner used in POLY-DEF, we specifically select a Ritz-based polynomial over a standard Chebyshev polynomial. This choice is twofold. First, constructing the Chebyshev polynomial requires an accurate estimate of $\lambda_{\min}$, which is computationally expensive for this highly ill-conditioned operator. Second, the Ritz-based polynomial dynamically adjusts to the specific discrete spectral distribution of the matrix, which empirically yields a considerably less dense accumulation of low modes in the resulting preconditioned operator $Ap(A)$.

Despite this polynomial selection, applying explicit LR-deflation to the preconditioned operator remains mandatory to safeguard the continuous contour quadrature. This deflation must target the critical Small Absolute Eigenmodes: Because the original matrix A possesses eigenvalues extremely close to the origin, the closed integration contour $\Gamma$ must cross the real axis at a very small positive value to safely bypass the branch cut. At these evaluation nodes, the mapped contour $\eta_{z}=zp(z)$ approaches zero. If the projected matrix $\tilde{H}_{m}$ retains un-deflated eigenmodes with small absolute values, the distance to singularity $|\eta_{z}-w|$ vanishes, causing a breakdown of the resolvent evaluation $(\eta_{z}I_{m}-T)^{-1}$. By deflating these well-resolved modes, we physically enforce a strictly bounded safe distance between the continuous evaluation path $\eta_{z}$ and the discrete spectrum of $\tilde{H}_{m}$

\begin{figure}[htbp]
    \centering
    \includegraphics[width=0.85\textwidth]{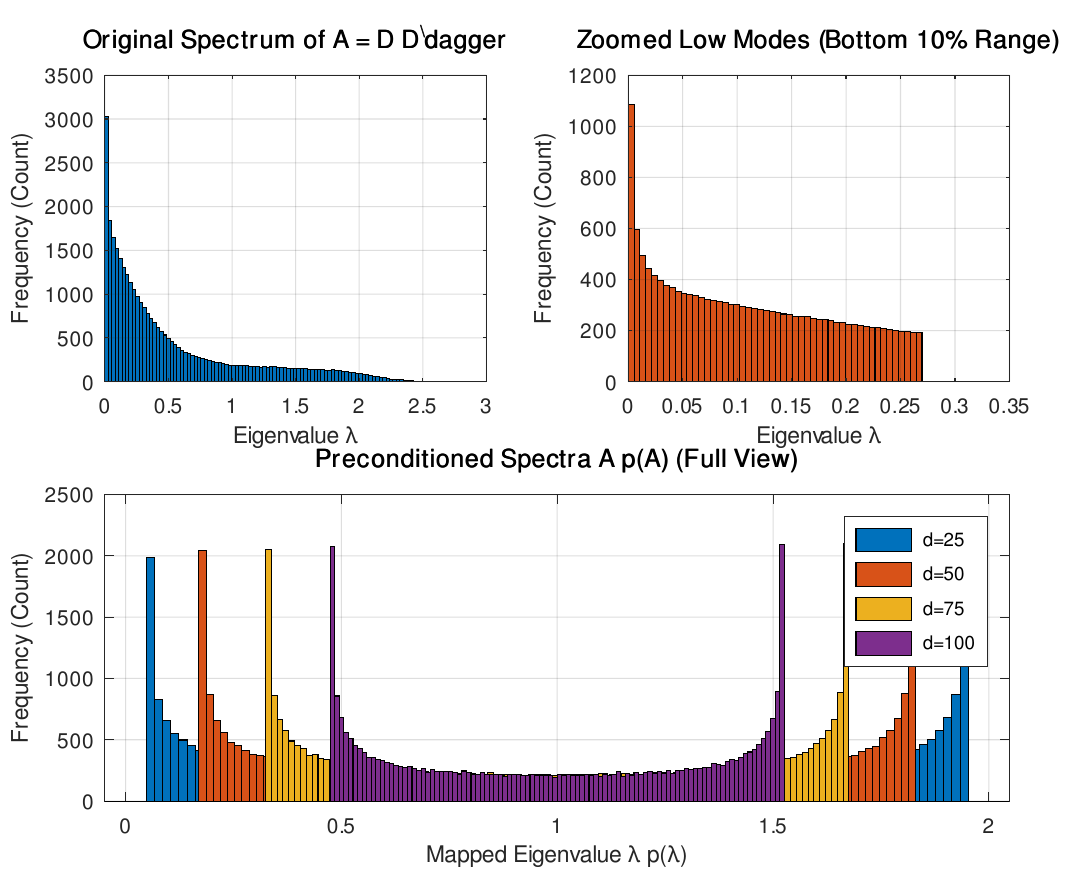}
    \caption{Spectral compression achieved by a standard Chebyshev polynomial for degrees $d \in \{25, 50, 75, 100\}$. The rigid bounds leave long tails and struggle to optimally cluster the low-frequency modes.}
    \label{fig:chebyshev_spectra}
\end{figure}

\paragraph{Spectral Clustering: The Advantage of Ritz over Chebyshev Polynomials} To evaluate the spectral compression achieved by these two polynomial preconditioners, we first examine the standard Chebyshev polynomial. Figure \ref{fig:chebyshev_spectra} illustrates the spectrum of the Chebyshev-preconditioned operator for degrees $d\in\{25,50,75,100\}$ \cite{ramirez23}. While the Chebyshev polynomial successfully clusters the bulk of the spectrum towards 1.0, its rigid formulation requires a priori expensive estimation of the spectral bounds $[\lambda_{\min},\lambda_{\max}]$. Furthermore, for extreme condition numbers, it leaves long, sparsely populated tails extending far from the target cluster, which limits its ability to optimally precondition the most critical low-frequency modes.

To improve low-mode clustering without relying on external bound estimations, we transition to Ritz-based polynomials \cite{reichel1991application,loe22}. The roots $\theta_{i}$ are extracted dynamically from the active Krylov subspace. To ensure orthogonality and prevent finite-precision basis degradation during the extraction, we employ a Double Modified Gram-Schmidt (DMGS) Arnoldi process.

\begin{figure}[htbp]
    \centering
    \begin{subfigure}[b]{0.48\textwidth}
        \includegraphics[width=\linewidth]{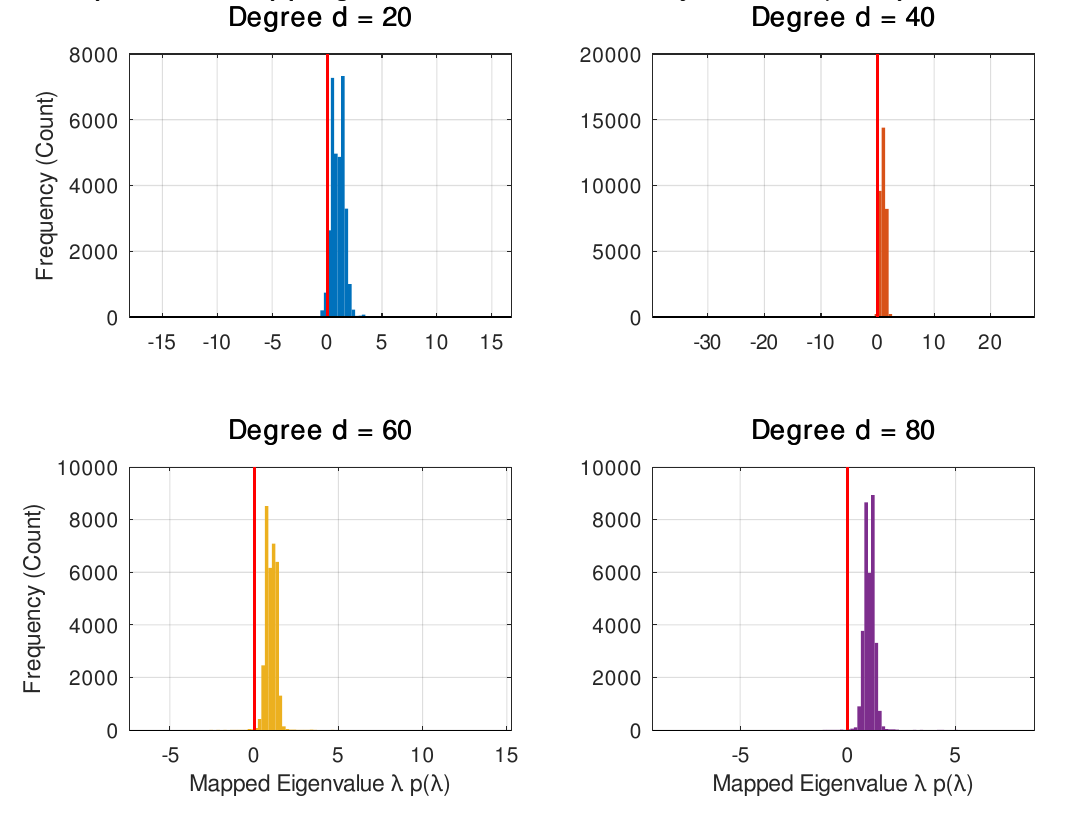}
        \caption{Non-Harmonic Ritz}
    \end{subfigure}
    \hfill
    \begin{subfigure}[b]{0.48\textwidth}
        \includegraphics[width=\linewidth]{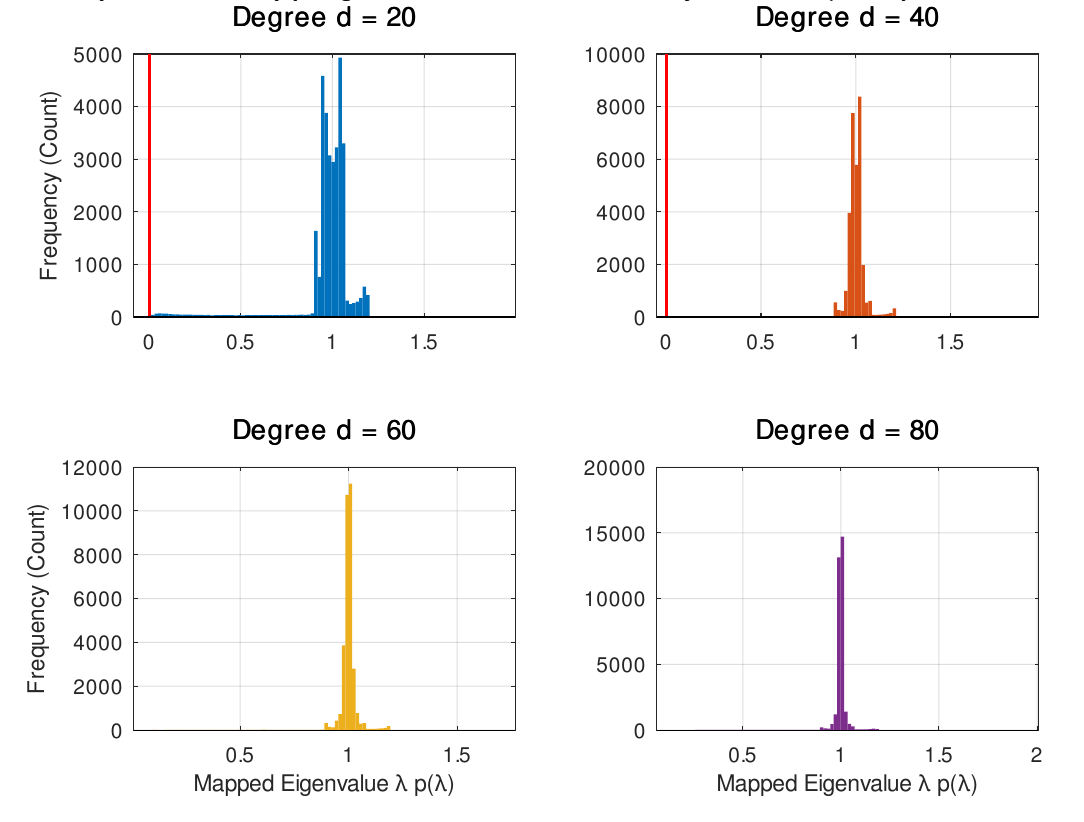}
        \caption{Harmonic Ritz}
    \end{subfigure}
    \caption{Comparison of spectral compression using non-harmonic (left) vs. harmonic (right) Ritz polynomials for varying degrees. Harmonic extraction natively targets near-origin eigenvalues, drastically improving global clustering toward 1.0.}
    \label{fig:ritz_comparison}
\end{figure}

Figure \ref{fig:ritz_comparison} directly compares the spectral compression of standard (non-harmonic) Ritz extraction against Harmonic Ritz extraction. The standard Ritz extraction (Figure \ref{fig:ritz_comparison}, left) intrinsically favors the largest eigenvalues, thereby starving the near-origin region of roots. This results in a highly dispersed preconditioned spectrum that fails to tightly cluster around 1.0. Conversely, the Harmonic Ritz extraction (Figure \ref{fig:ritz_comparison}, right) is designed to accurately locate interior and near-zero eigenvalues. By anchoring the polynomial roots precisely where the operator is most ill-conditioned, Harmonic extraction drastically improves the global clustering of the spectrum, pulling the bulk of the eigenvalues significantly closer to 1.0.

\begin{figure}[htbp]
    \centering
    \includegraphics[width=0.6\textwidth]{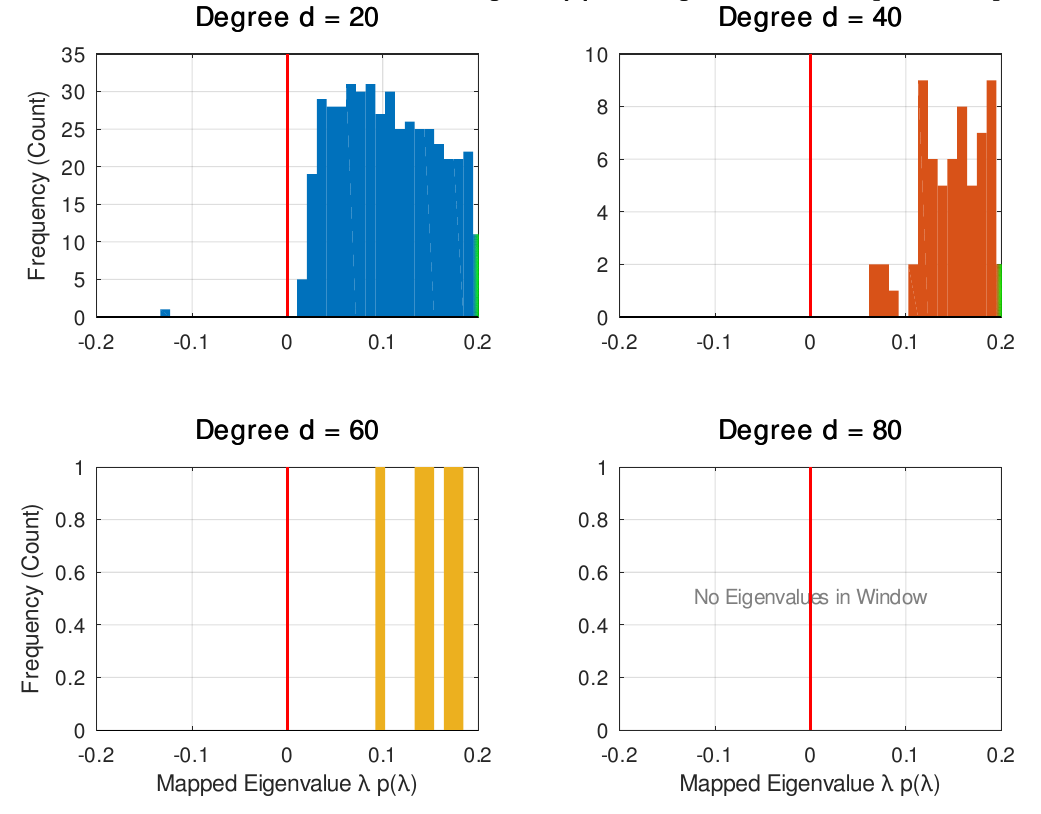}
    \caption{Zoomed-in view of the 30 smallest-real mapped eigenvalues for the Harmonic Ritz preconditioner. The remaining scattered low modes provide an ideal, well-resolved target for exact LR-deflation.}
    \label{fig:harmonic_zoomed}
\end{figure}

Figure \ref{fig:harmonic_zoomed} provides a high-resolution, zoomed-in view of the critical low-mode distribution (the 30 smallest-real mapped eigenvalues) for the Harmonic Ritz preconditioned operator. While Harmonic Ritz extraction excels at global clustering, this zoomed perspective reveals that a small number of eigenvalues inevitably remain scattered near 0.0. However, rather than representing a failure of the preconditioner, this specific spectral structure is highly advantageous. It isolates the most challenging modes, providing an ideal, well-resolved target for our exact LR-deflation mechanism (Section \ref{sec:lr_deflation}). By explicitly deflating these scattered low modes, the remaining spectrum is left clustered, safeguarding the continuous contour quadrature from singular poles and guaranteeing rapid, robust Krylov convergence.

\paragraph{Geometric Safety and Real-Axis Crossings} Because the base operator $A=DD^{\dagger}$ is Hermitian Positive Definite, the preconditioned operator $Ap(A)$ is also Hermitian, meaning its eigenvalues $\lambda_{p}$ are real. Consequently, the shifted preconditioned operator $\hat{A}_{z}=\eta_{z}I-Ap(A)$ can only become singular if the mapped contour scalar $\eta_{z}$ evaluates to a purely real number that coincides with an eigenvalue $\lambda_{p}$. Therefore, the continuous numerical quadrature is only threatened at the specific points where the complex mapped contour $\eta_{z}(\Gamma)$ crosses the real axis.

\begin{figure}[htbp]
    \centering
    \includegraphics[width=0.85\textwidth]{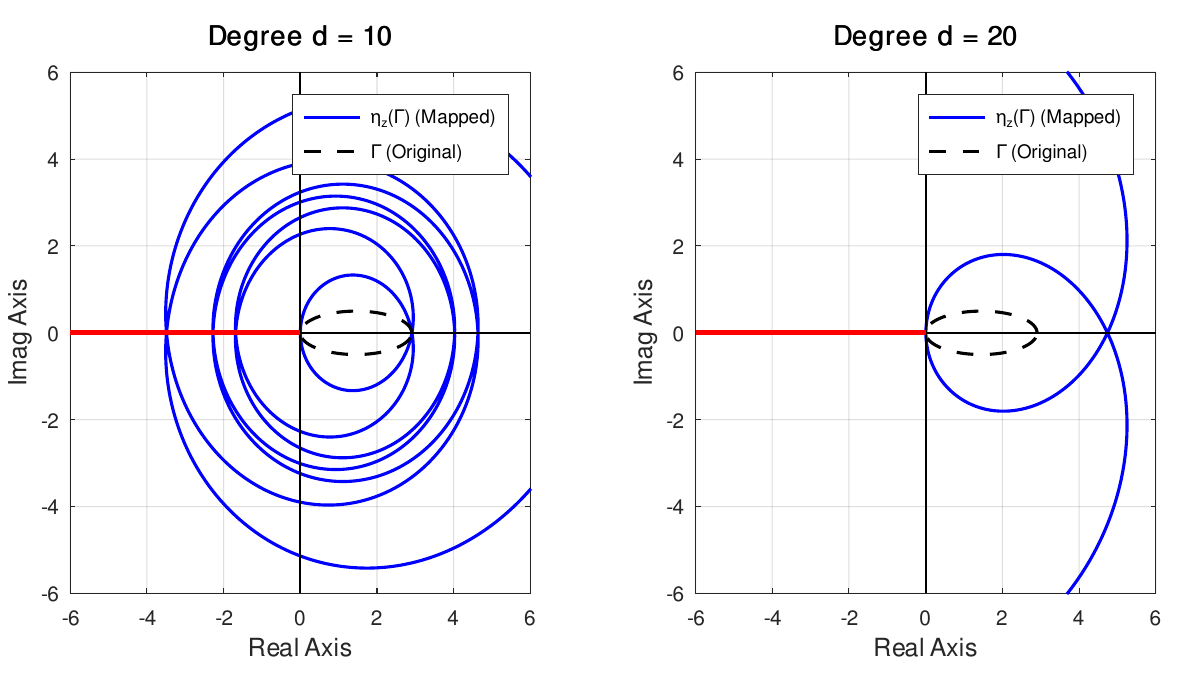}
    \caption{Mapping of the original elliptical contour $\Gamma$ (dashed black) to $\eta_{z}(\Gamma)$ (solid blue) for polynomial degrees $d=10$ (left) and $d=20$ (right). As d increases, the mapped contour aggressively expands outward, significantly increasing the safety margin for the numerical quadrature.}
    \label{fig:mapped_contours}
\end{figure}

To visualize this geometric interaction, Figure \ref{fig:mapped_contours} illustrates the mapping of the original closed elliptical contour $\Gamma$ through the polynomial shift parameter $\eta_{z}=zp(z)$. For a moderate polynomial degree such as $d=10$ (Figure \ref{fig:mapped_contours}, left), the mapped contour extends into the left-half of the complex plane, crossing the negative real axis in multiple locations. For a higher degree, $d=20$ on the right panel in that same figure, the risk of crossing the negative part of the real axis relatively close to the origin is much less. Higher degrees are therefore preferred in this regard.

In practice, there could be a clash between the mapped contour $\eta_{z}(\Gamma)$ and some eigenvalues of the preconditioned operator scattered towards the negative real axis. Although we do not see this in practice, and in our particular LQCD example, for $d \gtrsim 20$, we nevertheless deflate the smallest real (``sr" in Octave) to avoid such possible collisions. This guarantees that the active real spectrum of $Ap(A)$ never intersects the real-axis crossings of $\eta_{z}(\Gamma)$, securing a positive geometric safety margin $\min |\eta_{z}-\lambda_{p}|>0$ across the entire integration path.



\paragraph{Stable Evaluation via the Newton Basis} Crucially, executing the polynomial framework for the large degrees required here $(d\in[20,100])$ requires the deployment of the Newton basis formulation derived in Section \ref{sec:newton_basis}.

As established, evaluating either the mapped contour $(\eta_{z})$ or the global matrix assembly $(f_{m})$ using the standard monomial expansion $(\sum\gamma_{j}A^{j})$ triggers catastrophic numerical breakdown. The alternating coefficients grow exponentially with $d$, severely distorting the geometric contour and causing the intermediate vectors to swell by factors of $\mathcal{O}(10^{20})$.

By implementing POLY-DEF in the Newton root basis (naturally present in the Ritz-based polynomial), we eliminate the unstable coefficients $\gamma_{j}$ entirely for both the $\eta_{z}$ and the overall $f_{m}$ computation. Every operation—the contour mapping $\eta_{z}$, the dual-regime scalings $\tilde{\beta}_{j}(z)$, the LR-deflated scalar evaluation $N_{j}(\Lambda_{p})$, and the global full-space sequence $N_{j}(A)w_{j}$—is evaluated exactly via the bounded roots $\theta_{i}$ naturally extracted by the Harmonic Arnoldi process. This ensures numerical stability, protecting the framework against intermediate swell and finite-precision cancellation regardless of how large the polynomial degree $d$ becomes.

\paragraph{Algorithm Convergence and Deflation Efficiency}
To evaluate the performance, we track the relative convergence using the matrix-free stagnation proxy $\epsilon_m$ (evaluated solely on the Krylov component $w_{m,0,\mathrm{kry}}$ as derived in Section \ref{sec:stop_criterium_POLY_DEF}) against a target tolerance of $1.0 \times 10^{-10}$. For both configurations, i.e., with and without deflation, the preconditioner employs a Harmonic Ritz polynomial of degree $d=40$, stabilized via Leja ordering. 

\begin{figure}[htbp]
    \centering
    \begin{subfigure}[b]{0.48\textwidth}
        \includegraphics[width=\linewidth]{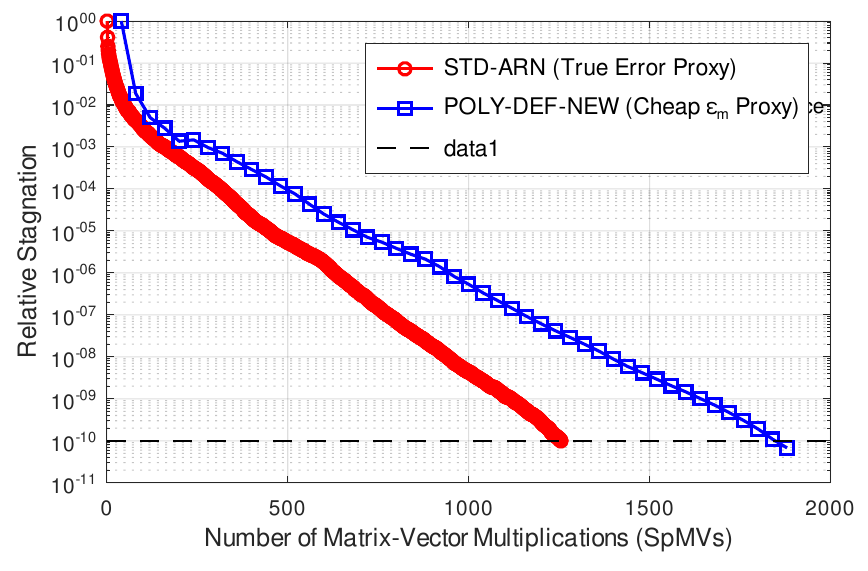}
        \caption{Without Deflation ($p=0$)}
    \end{subfigure}
    \hfill
    \begin{subfigure}[b]{0.48\textwidth}
        \includegraphics[width=\linewidth]{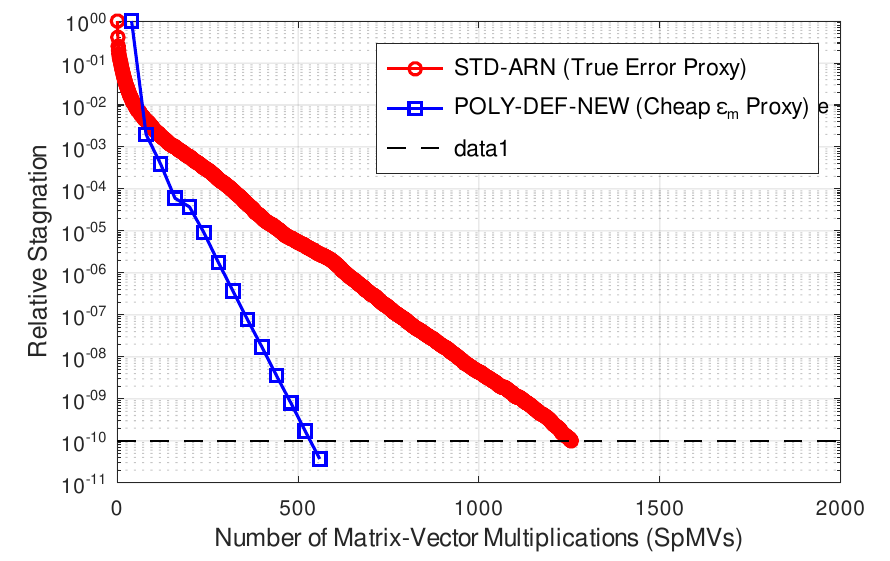}
        \caption{With Exact LR-Deflation ($p=400$)}
    \end{subfigure}
    \caption{Algorithmic convergence on the non-Hermitian twisted mass operator ($d=40$). While $p=0$ (left) requires a large $m \approx 1300$ baseline subspace to beat the polynomial in SpMVs, deflating 400 modes (right) achieves superior SpMV efficiency ($\approx 15$ Arnoldi steps) and allows a tenfold reduction in continuous quadrature nodes.}
    \label{fig:tm_convergence}
\end{figure}

Figure \ref{fig:tm_convergence} compares the convergence profiles in terms of total SpMV applications. In the absence of deflation ($p=0$, Figure \ref{fig:tm_convergence}a), the closed contour $\Gamma$ is forced to cross the real axis at a very small distance from the origin ($0.5 \times 10^{-4}$, with the upper crossing at $2.8$) to safely enclose the extreme low modes. Safely resolving these near-singular poles requires a high numerical quadrature resolution of $10^6$ integration points. From a SpMV-centric perspective, the baseline STD-ARN appears to outperform the preconditioned method here. However, this metric masks the severe memory and orthogonalization bottlenecks: STD-ARN requires an intractable Krylov subspace of dimension $m \approx 1300$, whereas POLY-DEF achieves convergence in merely $m \approx 47$ iterations ($\approx 1900$ SpMVs). Thus, the preconditioned method remains vastly superior in total wall-clock execution time due to the elimination of $\mathcal{O}(N m^2)$ orthogonalization costs, even once one switches to a double-pass Lanczos.


The true advantage of the framework is evident when exact LR-deflation is deployed ($p=400$, Figure \ref{fig:tm_convergence}b). We extract the 400 smallest real (\texttt{'sr'}) eigenmodes of the preconditioned operator $AM_{0}^{-1}$ to a relative tolerance of $1.0{\cdot}10^{-10}$. By removing these critical modes from the active Krylov projection, the contour's lower real-axis crossing is safely shifted outward to $1.0{\cdot}10^{-4}$. This widened geometric safety margin permits a tenfold reduction in the continuous quadrature requirement, down to $10^{5}$ nodes. As a result, the deflated POLY-DEF method exhibits significant acceleration, reaching the desired tolerance in fewer than 600 SpMVs ($\approx15$ Krylov iterations). This outperforms the standard Arnoldi baseline in raw operator applications, dot products and dense subspace algebra.

\section{Conclusions and Future Work}

The presented framework offers two robust paths for preconditioning $f(A)b$. The rational shift-and-invert approach yields an exact closed-form extraction that avoids contour integration, highly favorable when sparse linear system solves are permissible via multigrid. By rigorously shielding the subspace projection via Double Modified Gram-Schmidt reorthogonalization, we eliminate finite-precision phantom poles and resolve the fundamental spectral trade-off. Conversely, the polynomial framework provides a matrix-free approach. By decoupling the preconditioner via shift-invariance, incorporating stable Schur dual-regime quadrature, and insulating against contour singularities via LR-deflation, it offers a highly parallelizable algorithm specifically engineered to evaluate $f(A)b$ with optimal SpMV efficiency.

It is worth highlighting a critical algorithmic advantage of the proposed polynomial framework over previous approaches \cite{frommer2024polynomial}. In prior methodologies, evaluating the matrix inverse square root required algebraically extracting the polynomial out of the function, which introduced a problematic sign choice when evaluating the inverse square root of the squared polynomial. The contour integral framework presented here bypasses this issue entirely. Because the polynomial is built directly on $A = DD^\dagger$ and the continuous formulation natively incorporates the shifted preconditioner, we no longer need to formally take the inverse square root of the square of the polynomial itself, ensuring a robust, unambiguous evaluation.

As future work, we plan to deploy both branches of the framework on realistic, fine-grid Lattice QCD (LQCD) problems. Specifically, we will apply the matrix-free POLY-DEF method to large-scale LQCD configurations, and concurrently explore the deployment of SAI-MG on the exact same examples by leveraging DD$\alpha$AMG as the internal sparse linear solver. 

When using POLY-DEF, although the extraction of low-modes to be deflated is done on the preconditioned operator, leading to various computational advantages in terms of eigensolving execution time, one might need to compute too many such low-modes, leading to a method that can simply not win even over standard Arnoldi. For such cases, a possible direction of future research is the integration of deflation/recycling strategies \cite{burke2022krylov} within the framework presented here.

Furthermore, we intend to explore a broader class of preconditioners by relaxing the operator equivalence condition defined in eq.\ \ref{eq:krylov_subspace_shift_invariance}. This relaxation could open the door to novel, highly efficient preconditioning structures. Finally, a very interesting theoretical direction is the investigation of preconditioners that fundamentally break shift-invariance altogether. Such an approach would preclude the use of the standard shifted approximations used here, requiring the use of an entirely different Krylov projection strategy, and representing a compelling frontier for accelerating the evaluation of general matrix functions.

\bibliographystyle{siamplain}
\bibliography{references}

\end{document}